\documentstyle[12pt]{article}
\begin{document}


\def\c{{\bf C}}
\def\n{{\bf N}}
\def\z{{\bf Z}}
\def\u{{\bf U}}
\def\f{{\bf F}}
\def\t{{\bf T}}
\def\s{{\bf S}}
\def\o{{\bf O}}
\def\w{{\bf W}}
\def\g{{\bf G}}
\def\r{{\bf R}}
\def\ll{{\bf L}}

\def\ab{\overline{a}}
\def\bb{\overline{b}}
\def\cb{\overline{c}}
\def\db{\overline{d}}
\def\zb{\overline{z}}
\def\ub{\overline{u}}
\def\vb{\overline{v}}
\def\fb{\overline{F}}
\def\rb{\overline{r}}
\def\xb{\overline{x}}
\def\Hb{\overline{H}}

\def\no{\mid\!\mid}
\def\lvs{\vskip 1mm}
\def\qed{\Box}
\def\l{{\cal L}}
\hyphenation{Wo-ro-no-wi-cz}
\hyphenation{a-me-na-ble}

\newtheorem{defi}{Definition}[section]
\newtheorem{prop}{Proposition}[section]
\newtheorem{theo}{Theorem}[section]
\newtheorem{lemm}{Lemma}[section]
\newtheorem{coro}{Corollary}[section]

\noindent {\Large\bf Representations of compact quantum groups
and subfactors}

\bigskip\noindent {\large Teodor Banica}

\bigskip\noindent {\bf Abstract:} We associate Popa
systems (= standard invariants of subfactors) to the finite dimensional representations
of compact quantum groups. We characterise the systems arising
in this way: these are the ones which can be ``represented'' on
finite dimensional Hilbert spaces. This is proved by a universal
construction. We explicitely compute (in terms of some free products) the operation of going from representations
of compact quantum groups to Popa systems and then back via the
universal construction. We prove a Kesten type result for the co-amenability
of compact quantum groups, which allows us to compare it with
the amenability of subfactors.


\section*{Introduction}

In this paper we investigate the relationship between representations
of compact quantum groups and standard invariants of subfactors. We
construct and study three operations connecting these objects (see the abstract) and we discuss the notion of amenability. In a certain sense,
these are the only universal operations connecting representations of
compact quantum groups and Popa systems (see the concluding remarks in the end of this section).

A few words on the terminology. The paper is written using the
formalism of Woronowicz algebras, which are the Hopf
$\c^*$-algebras which  correspond to both notions of ``algebras of
continuous functions on compact quantum groups'' and ``$\c^*$-algebras
of discrete quantum groups''. Part of the results are better understood in
terms of compact quantum groups, and the other part in terms of
discrete quantum groups, so we have written this introductory
section by using the more suggestive formalism of quantum
groups associated with Woronowicz algebras.

It has to be mentioned that the locally compact quantum
groups (\cite{bs}) are known to be related to the irreducible depth 2
subfactors by a crossed product construction (cf. Ocneanu's theorem and its
extensions; see \cite{en} and the references therein). The point of
view in this paper is different - in fact we relate the ``theory of a
single representation of compact quantum groups'' with the ``theory of
the standard invariants of subfactors''. 

We recall that if $N\subset M$ is a finite index inclusion of $II_1$
factors, its {\em lattice of higher relative commutants} is by
definition the lattice $\{ M_i^\prime\cap M_j\}_{0\leq i\leq
  j<\infty}$, where $\{ M_i\}_{i\geq 0}$ denotes the sequence of
algebras in the Jones tower. Here by {\em lattice of $\c^*$-algebras}
(or just {\em lattice}) we will mean a collection of $\c^*$-algebras $\{ A_{ij}\}_{0\leq i\leq j<\infty}$,
together with given inclusions between them, as follows:
$$\matrix{
A_{00} &\subset & A_{01} & \subset & A_{02}&
\subset & A_{03}& \subset\
\cdots\cr &\ &\ \cup &\ &\cup &\ &\cup \cr &\ &\ A_{11} & \subset &
A_{12}&
\subset & A_{13}& \subset\ \cdots\cr &\ &\
&\ &\cup &\ &\cup \cr &\ &\ \ & \ & A_{22} & \subset & A_{23}&
\subset\ \cdots\cr &\ &\ &\ &\ &\ &\cup \cr
&\ &\ &\ &\ &\ &\cdots &\ \ \ \ \ \cdots\cr
}$$
The lattice of higher relative commutants, together with its additional structure coming from traces and Jones projections is
called the {\em standard invariant} of the subfactor. It is a
complete invariant under amenability assumptions
(see \cite{p25}). When regarded as an abstract object satisfying the axioms
in \cite{p3}, it will be called {\em Popa system}.

We start from
the general principle that if $\pi$ is a ``representation of a quantum
group'' then the following lattice $L(\pi )$ of algebras 
$$\matrix{
\c &\subset & End(\pi ) & \subset & End(\pi\otimes \hat{\pi})&
\subset & End(\pi\otimes \hat{\pi}\otimes \pi )& \subset\
\cdots\cr &\ &\ \cup &\ &\cup &\ &\cup \cr &\ &\ \c & \subset &
End(\hat{\pi})&
\subset & End(\hat{\pi}\otimes \pi )& \subset\ \cdots\cr &\ &\
&\ &\cup &\ &\cup \cr &\ &\ \ & \ & \c & \subset & End(\pi )&
\subset\ \cdots\cr &\ &\ &\ &\ &\ &\cup \cr
&\ &\ &\ &\ &\ &\cdots &\ \ \ \ \ \cdots\cr
}$$
should be the standard invariant of a subfactor. A heuristic
explanation for this fact is provided for instance by Ocneanu's
bimodule picture of the standard invariant: the standard invariant of
any subfactor $N\subset M$ is ``of the above form'', with $\pi$
replaced by the bimodule ${\ }_N\!L^2(M)_M$ (see e.g. \cite{bi}).

This general fact was intensively investigated for the quantum
groups at roots of unity (which of course are not compact
quantum groups), and a whole series of interesting subfactors
was constructed in this way, see
\cite{we1},\cite{we2},\cite{we3},\cite{xu} and also
\cite{jw},\cite{was2},\cite{to}.

The above general statement holds also when $\pi$ is a representation
of a deformation $G_q$ of a compact group (with $q>0$) \cite{xu}. Also
Wassermann's subfactors associated to representations of compact
groups (\cite{was}) and the locally trivial subfactors of Jones
and Popa (see \cite{p1}) have higher relative commutants of the form
$L(\pi )$, with $\pi$ a representation of a compact group $G$
(resp. of a ``dual'' $\widehat{\Gamma}$ of a discrete group). These objects $G_q$
($q>0$), $G$ and $\widehat{\Gamma}$ are compact quantum groups
in the sense of Woronowicz \cite{w1},\cite{w2},\cite{w3} (for
$G_q$, see
\cite{ro}), and in fact the following general result holds.

\lvs
\noindent {\bf Theorem A.} {\em If $\pi$ is a finite dimensional
unitary representation of a compact quantum group $G$, then the
lattice $L(\pi )$ is a Popa system.}

\lvs
That is, $L(\pi )$ has to be the lattice of higher relative
commutants of a subfactor (by \cite{p3}). The verification of
Popa's axioms is very easy when the Haar measure $h\in C(G)^*$
is a trace; in the general case one has to perturb everything
by using the characters $f_z\in {\cal R}(G)^*$ describing the
modular theory of the Haar measure \cite{w1}. The resulting
index is the square of the quantum dimension of $\pi$.

There are several restrictions on the Popa systems arising in
this way, and the most obvious one is that if $H$ is the
Hilbert space where $\pi$ acts, then $\hat{\pi}$ acts on
$\Hb$, $\pi\otimes\hat{\pi}$ acts on $H\otimes \Hb$ etc., so
$L(\pi )$ has to be a sublattice of the following
lattice $L(H)$:
$$\matrix{
\c &\subset & \l (H) & \subset & \l (H\otimes \Hb )& \subset & \l (H\otimes \Hb\otimes H)& \subset\ \cdots\cr
&\ &\ \cup &\ &\cup &\ &\cup \cr
&\ &\ \c & \subset & \l (\Hb )& \subset & \l (\Hb\otimes H)& \subset\ \cdots\cr
&\ &\ &\ &\cup &\ &\cup \cr
&\ &\ \ & \ & \c & \subset & \l (H)& \subset\ \cdots\cr
&\ &\ &\ &\ &\ &\cup \cr
&\ &\ &\ &\ &\ &\cdots &\ \ \ \ \ \cdots\cr
}$$
{\bf Theorem B.} {\em For any Popa
system $(A_{ij})_{0\leq i\leq j<\infty}$ the following are
equivalent:

\noindent (i) there exists a compact quantum group $G$ and
a finite dimensional representation $\pi$ of $G$ such that
$(A_{ij})_{0\leq i\leq j<\infty}$ is equal to $L(\pi )$.

\noindent (ii) there exists a finite dimensional Hilbert
space $H$ such that $(A_{ij})_{0\leq i\leq j<\infty}$ is a
sublattice of $L(H)$, and such that the Jones projection $e_2\in
A_{02}$ corresponds in this way to a rank one projection in $\l
(H\otimes \Hb)$.}

\lvs 
The part {\it (i)} $\Longrightarrow$ {\it (ii)} is trivial, and the
proof of the converse is based on the
following observation. Associated to any compact quantum group
$G$ is its monoidal category of representations $rep(G)$,
together with the forgetful functor $F_G:rep(G)\rightarrow
vect_\c$. The Tannakian duality \cite{w2} says that the
compact quantum groups are in one-to-one correspondence with
the pairs (monoidal category, monoidal functor) satisfying
a certain list of axioms. Now if $\pi$ is a representation of
$G$, the algebras in $L(\pi )$ could be thought of as a ``piece of
$rep(G)$''; and the embedding of $L(\pi )$ into $L(H)$ could be
thought of as the ``corresponding piece of $F_G$''. Thus what we
have to do for proving {\it (ii)} $\Longrightarrow$ {\it (i)}
is to construct a monoidal category and a monoidal functor
when knowing ``pieces'' of them; and this kind of problem
(see also \cite{w2},\cite{kw},\cite{b0},\cite{b1}) is
well-known to be usually a combinatorial one. 

The construction is not unique, and in fact we find the
``universal'' pair $(G,\pi )$ such that $(A_{ij})_{0\leq i\leq
j<\infty}$ is equal to $L(\pi )$ (as sublattices of $L(H)$). The
following question arises naturally: if $(A_{ij})_{0\leq i\leq
j<\infty}$ itself is the Popa system associated to a pair
$(G,\pi )$, who is the universal pair that we construct ? 

\lvs
\noindent {\bf Theorem C.} {\em If $\pi$ is a representation
of $G$ on $H$ and if $(\tilde{G},\tilde{\pi})$ is the universal
pair satisfying $L(\tilde{\pi})=L(\pi )$ as sublattices of
$L(H)$ then there exists an (explicit) embedding
$\widehat{\tilde{G}}\hookrightarrow \z *\widehat{\pi (G)}$.}

\lvs
Here $*$ is the free product of discrete quantum groups
\cite{wa} (see section 5 for the rigorous statement, in
terms of Hopf algebras). The proof uses an isomorphism
criterion which relies on a result of free probability theory
of Nica and Speicher \cite{ns} and on the following idea from
\cite{b1}: the dimensions of the linear spaces $Hom(r,p)$
with $r,p=$ tensor products between $\pi$ and $\hat{\pi}$ are
exactly the $*$-moments of the character $\chi (\pi )\in C(G)$
with respect to the Haar measure $h\in C(G)^*$.

The above considerations give also a more conceptual proof for
the result $A_u(F)_{red}\hookrightarrow\c^*(\z )*_{red}A_o(F)$
from \cite{b1}. In fact there are also other results from
\cite{b1} on $A_u(F)$ which extend to the algebras of the form
$C(\tilde{G})$.

The last result is about amenability. For locally compact
quantum groups the basic results on amenability were
established by Blanchard \cite{bl}, and in the discrete case a
Kesten type result could be deduced from his work (this was
explained to us by Skandalis). This result has
several applications (see section 6), one of them being:

\lvs
\noindent {\bf Theorem D.} {\em If $\pi$ is a
representation of $G$ then $L(\pi )$ is amenable in the sense
of \cite{p25} if and only if $K=(\pi\otimes\hat{\pi })(G)$ has
the following properties: the dual $\widehat{K}$ is amenable,
and the Haar measure $h\in C(K)^*$ is a trace.}

\lvs 
This generalises some known results, and also shows that the
index of an amenable lattice of the form
$L(\pi )$ has to be the square of an integer. 

\lvs 
\noindent {\bf Concluding remarks.} The theorems A,B,C
could be interpreted in the following way. Fix $H$ and let $X$
be the category of pairs (compact quantum group, representation on
$H$) and $Y$ be the set of pairs (Popa system, embedding into
$L(H)$). Then A,B,C (and their proofs) give respectively: a
surjection $X\rightarrow Y$; the universal section $Y\rightarrow
X$; a description of the corresponding projection $X\rightarrow
X$. It follows that the elements of $Y$ are in one-to-one
correspondence with the elements in the image of the
projection. 

One remaining problem on the relationship between $X$ and Popa
systems is to find the relationship between $Y$ and Popa
systems.

One can show that the universal construction in section 4
makes sense for any Popa system $(A_{ij})_{0\leq i\leq
j<\infty}$, by giving abstract meanings to the operators
$i$ and $p$ (e.g. small semicircles when working with planar
diagrams \cite{jo}). If ${\cal R}_{(A_{ij})_{0\leq i\leq
j<\infty} }$ is the monoidal category constructed in this way,
then the above question of existence, uniqueness, classification
of embeddings of $(A_{ij})_{0\leq i\leq j<\infty}$ into
$L(H)$'s is equivalent to the question of existence, uniqueness,
classification of fiber functors on
${\cal R}_{(A_{ij})_{0\leq i\leq j<\infty} }$. While this
can be done in certain cases, in general this is related to
an unsolved problem, namely the generalisation of the
theorems of Doplicher-Roberts \cite{dr} and Deligne \cite{de}. 

\lvs 
The paper is organised as follows. In the first section we
recall Woronowicz' formalism, we discuss the notion of
duality for unitary representations, and we give a list of
relevant examples. In the second section we prove the theorem A
and we give some examples. In the third and the fourth
section we prove the theorem B. In the fifth and the sixth
section we prove the theorems C, respectively D.

\lvs
Part of this work was done during my stay at UCLA and
the University of Gen\`eve; I would like to thank these
institutions and Sorin Popa for their warm hospitality. I am
also grateful to Dietmar Bisch and Georges Skandalis for
useful discussions. 

\section{Duality for corepresentations}

The Woronowicz algebras are the Hopf $\c^*$-algebras which 
correspond (by \cite{w1},\cite{w2},\cite{w3},\cite{bs} etc.) to
both notions of ``algebras of continuous functions on compact
quantum groups'' and ``$\c^*$-algebras of discrete quantum
groups''. They can be defined as being the bisimplifiable
unital Hopf $\c^*$-algebras, see \cite{w3}. One alternative
definition, which will be used in this paper (see the note
below), is as follows. 

\lvs 
Consider pairs $(A,u)$ consisting of a unital $\c^*$-algebra
$A$ and a unitary matrix $u\in M_n(A)$ subject to the
following  conditions (\cite{w1}, definition 1.1):

{\it (w1)} the coefficients of $u$ generate
in $A$ a dense $*$-subalgebra, called $A_s$.

{\it (w2)} there exists a $\c^*$-morphism $\delta
:A\rightarrow A\otimes_{min}A$ sending 
$u_{ij}\mapsto\sum u_{ik}\otimes u_{kj}$ for any $1\leq
i,j\leq n$.

{\it (w3)} there exists a linear antimultiplicative
map $\kappa :A_s\rightarrow A_s$ such that 
$\kappa (\kappa (a^*)^*)=a$ for all $a\in A_s$ and such that 
$(id\otimes\kappa )u=u^{-1}$. 

If $(A,u)$ and $(B,v)$ are as above, one can define the
morphisms from $(A,u)$ to $(B,v)$ to be the $\c^*$-morphisms
$f:A\rightarrow B$ satisfying $\delta_Bf=(f\otimes
f)\delta_A$. The {\em Woronowicz algebras} are by
definition the inductive limits in the corresponding category of
pairs $(A,u)$ satisfying {\it (w1-3)}.

\lvs 
If $(A,u)$ satisfies {\it (w1-3)} then the dense subalgebra
$A_s$ is an involutive Hopf
$\c$-algebra with comultiplication $\delta$, antipode
$\kappa$ and counit defined by $\varepsilon
:u_{ij}\mapsto\delta_{ij}$. Recall that the finite dimensional
corepresentations of $A_s$ are the elements 
$r\in {\cal L}(V)\otimes A_s$ satisfying $(id\otimes\delta
)r=r_{12}r_{13}$, where $V$ is a finite dimensional complex
vector space (for instance the matrix $u$ is a unitary
corepresentation of $A_s$ on $\c^n$, called the fundamental
corepresentation of $(A,u)$). It was proved in \cite{w1} that
every non degenerate such corepresentation is equivalent to a
unitary corepresentation, and is completely reducible. The
Hopf algebra $A_s$ has then a (cosemisimple) decomposition
$$A_s=\bigoplus_{r\in Irr(A)} A(r)$$ 
where $Irr(A)$ is the set of equivalence classes of finite dimensional irreducible corepresentations of $A_s$, and 
for every (class of) corepresentation $r\in \l (V)\otimes A_s$ its linear space of coefficients is 
$A(r):=\{ (\varphi\otimes id)r\mid \varphi\in {\cal L}(V)^*\}$.

These results extend to all Woronowicz algebras (see
\cite{w3},\cite{bs},\cite{wa}). 

It is easy to see that the
finite dimensional corepresentations of $C(G)$ correspond to
the finite dimensional representations of $G$ (for $G=$ compact
group) and that the finite dimensional irreducible
corepresentations of
$\c^*(\Gamma )$ correspond to the elements of $\Gamma$ (for
$\Gamma =$ discrete group). 

\lvs 
\noindent {\bf Note.} If $u$ is a finite dimensional unitary
corepresentation of a Woronowicz algebra $B$, and if $A$ is the
$\c^*$-subalgebra of $B$ generated by the coefficients of $u$
then $(A,u)$ satisfies the conditions {\it (w1-3)}. As we will
be interested in the theory of a single corepresentation of a
Woronowicz algebra, all the Woronowicz algebras $A$ we consider
in this paper may be supposed to satisfy the conditions {\it
(w1-3)} for some matrix
$u\in M_n(A)$, i.e. may be supposed to be as in \cite{w1}.
Note that these Woronowicz algebras are exactly the ones whose
associated compact quantum groups are ``compact matrix quantum
groups'', or, equivalently, whose associated discrete quantum
groups are ``of finite type''.
\lvs 
One of the main results of \cite{w1} was the construction of a canonical family of characters 
$\{ f_z\}_{z\in\c}$ of $A_s$, which describe the modular
theory of the Haar measure $h\in A^*$. The $f_z$'s have the
following properties (where $*$ denotes the convolutions over
the Hopf algebra $A_s$):
\lvs 
{\it (f1)} $f_z*f_{z^\prime}=f_{z+z^\prime}$, 
$\forall\, z,z^\prime\in\c$ and $f_0=\varepsilon$ (the counit
of $A_s$). 

{\it (f2)} $f_z\kappa (a)=f_{-z}(a)$, $f_z(a^*)=\overline{f_{-\overline{z}}(a)}$, $\forall\, a\in A_s$ and $z\in\c$.

{\it (f3)} $\kappa^2(a)=f_{-1}*a*f_1$, $\forall\, a\in A_s$.

{\it (f4)} $h(ab)=h(bf_1*a*f_1)$, $\forall\, a,b\in A_s$.
\lvs 
If $v\in \l (H)\otimes A_s$ is a finite dimensional unitary
corepresentation of $A_s$ then the restriction of $f_z$ to the
space of coefficients of $v$ can be computed by using the
formula $(id\otimes f_z)v=Q_v^{2z}$, where
$$Q_v:=(id\otimes f_{\frac{1}{2}})v\in \l (H)$$
(this follows from {\it (f1)}). In lemma 1.1 we list the properties of
the operator $Q_v$, and in lemma 1.2 we give a characterisation of $Q_v$.  
\lvs 
\noindent {\bf Notations.} If $H$ is a finite dimensional Hilbert
space, $\Hb$ is its complex conjugate. We will use the
identifications $\overline{\Hb}=H$ and $\overline{H\otimes
K}=\overline{K}\otimes\Hb$ (but not $H\simeq \Hb$, even if
$H=\c^n$). The transposition $\l (H)\rightarrow
\l (\Hb )$ will be denoted by $a\mapsto t(a)$ or $a\mapsto
a^t$.
 
\begin{lemm}
If $v\in\l (H)\otimes A_s$ is a finite dimensional 
unitary corepresentation of $A_s$ then $Q_v$ has the following
properties:

(i) $Q_v>0$.

(ii) $Q_v\in End(v)^\prime$. 

(iii) $Tr(Q_v^2.)=Tr (Q_v^{-2}.)$ on $End(v)$.

(iv) $Q_v^2\in Hom(v,(id\otimes\kappa^2)v)$.

(v) $Q_v^t\,\vb\, (Q_v^t)^{-1}$ is unitary.
\end{lemm} 

{\em Proof.} We prove first (v). We have 
$$Q_v^t\,\vb\, (Q_v^t)^{-1}=[(t\otimes f_{\frac{1}{2}})v][(t\otimes\kappa )v][(t\otimes f_{-\frac{1}{2}})v]=(t\otimes j)v$$
where $j:A_s\rightarrow A_s$ is the linear map 
$x\mapsto f_{\frac{1}{2}}*\kappa (x)*f_{-\frac{1}{2}}$. 
It's easy to see using {\it (f1,f2,f3)} that $j$ is 
antimultiplicative, and that it commutes with the involution of
$A_s$. Thus $t\otimes j:\l (H)\otimes A_s\rightarrow \l
(\overline{H})\otimes A_s$ is an antimorphism of
$*$-algebras, so it maps unitaries to unitaries and this proves
(v). 

If $v$ is irreducible then (ii) is trivial and (i,iii,iv)
follow from  {\it (f1-f4)}, see section 5 in
\cite{w1}. The results (i-iv) can be extended to the
general case by using cosemisimplicity: by \cite{w1} we may
suppose that $v=\sum m_iv_i$ with $m_i\in\n$ and $v_i\in \l
(H_i)\otimes A_s$  irreducible and disjoint. As
$m_iv_i=v_i\otimes (m_i1)$ acts on $\l (H_i)\otimes \l
(\c^{m_i})$  for every $i$ we may suppose that $H=\oplus_i
(H_i\otimes \c^{m_i})$. With respect to this decomposition we
have
$$End(v)=\bigoplus_i id_{H_i}\otimes \l (\c^{m_i})$$ 
and $Q_v=\sum_i Q_{v_i}\otimes id_{\c^{m_i}}$, and this 
proves (i,ii,iv). Also if $x\in End(v)$ then $x$ is of
the form $\sum_i id_{H_i}\otimes x_i$, so 
$$Tr(xQ_v^2)=\sum_i Tr(Q_{v_i}^2)Tr(x_i)=\sum_i Tr(Q_{v_i}^{-2})Tr(x_i)=Tr(xQ_v^{-2})$$
and this proves (iii). $\qed$

\begin{lemm}
If $v\in\l (H)\otimes A_s$ is a finite dimensional unitary corepresentation of $A_s$ and $Q\in \l (H)$ is a positive operator 
such that $Tr(Q^2.)=Tr(Q^{-2}.)$ on $End(v)$, then the following are equivalent:

(i) $Q=Q_v$.

(ii) $Q^t\,\vb\, (Q^t)^{-1}$ is unitary.

(iii) $Q^2\in Hom(v,(id\otimes\kappa^2)v)$.

(iv) $QQ_v^{-1}\in End(v)$.
\end{lemm}

{\em Proof.} {\it (i} $\Longrightarrow$ {\it ii)} was part of the preceding lemma. 

{\it (ii} $\Longrightarrow$ {\it iii)} Using the formula $(id\otimes\kappa )r=r^*$ for the unitary corepresentations 
$r:=v$ and $r:=Q^t\, \vb\, (Q^t)^{-1}$ 
we get $(id\otimes\kappa )v=v^*$ and $(id\otimes\kappa )\vb =(Q^t)^{-2}\, v^t\, (Q^t)^2$, 
so $(id\otimes\kappa^2 )v=Q^2vQ^{-2}$. 

{\it (iii} $\Longrightarrow$ {\it iv)} We have $Q^2,Q_v^2\in Hom(v,(id\otimes\kappa^2)v)$, so $Q^2=Q_v^2x$ for some $x\in End(v)$. As 
$Q_v$ commutes with $x$, the result follows by taking square roots.

{\it (iv} $\Longrightarrow$ {\it i)} By using a decomposition $v=\sum v_i\otimes (m_i1)$ with $v_i\in \l (H_i)\otimes A_s$ 
as in the proof of the preceding lemma, we may write $QQ_v^{-1}=\sum_i id_{H_i}\otimes x_i$ with $x_i\in \l (\c^{m_i})$. Thus 
$Q=\sum Q_{v_i}\otimes x_i$, so for every $i$ the restriction to $End(v_i\otimes (m_i1))=\l (\c^{m_i})$ of the linear form 
$Tr(Q^2.)$ (resp. $Tr(Q^{-2}.)$) is $Tr(Q_{v_i}^2)Tr(x_i^2.)$ (resp. $Tr(Q_{v_i}^{-2})Tr(x_i^{-2}.)$). As 
$Tr(Q_{v_i}^2)=Tr(Q_{v_i}^{-2})$, we have $Tr(x_i^2.)=Tr(x_i^{-2}.)$ on 
$\l (\c^{m_i})$, so $x_i^2=x_i^{-2}$ for every $i$. But all $x_i$'s are positive, so they have to be the identities, and 
$Q=Q_v$. $\qed$
\lvs
\noindent {\bf Definition 1.1} (cf. lemma 1.1) 
If $v\in \l (H)\otimes A_s$ is a finite dimensional unitary corepresentation of $A_s$, the positive number 
$$d(v):=Tr(Q_v^2)=Tr(Q_v^{-2})=(Tr\otimes f_1)v=f_1\chi (v)$$
is called the quantum dimension of $v$. The linear form 
$$\tau_v:x\mapsto d(v)^{-1}Tr(Q_v^2x)=d(v)^{-1}Tr(Q_v^{-2}x)$$
is called the canonical trace on $End(v)$. The unitary corepresentation 
$$\hat{v}:=(Q_v)^t\,\vb\, (Q_v^t)^{-1}=
(t\otimes f_{\frac{1}{2}}*\kappa (.)*f_{-\frac{1}{2}})v$$
is called the canonical dual of $v$.
\lvs
Note that $\tau_v$ is a faithful positive unital trace on the 
$\c^*$-algebra $End(v)$. By the Cauchy-Schwartz inequality
applied to the eigenvalues of $Q_v^2$ and of $Q_v^{-2}$ we see
that $d(v)\geq dim(H)$, with equality iff $Q_v=id_H$. 

We will give now a list of examples of Woronowicz algebras,
with computation of the objects in the above definition,
and/or of the fusion semiring. We recall that if $A$ is a
Woronowicz algebra, the fusion semiring $R^+(A)$ is the
semiring whose elements are equivalence classes of finite
dimensional corepresentations of $A$, and whose operations are
the sum and tensor product of (classes of) corepresentations.
Its computability is of course the most important requirement
for having the computability of the (fusion algebra of) Popa
systems associated to $A$. 

The examples below illustrate most of the phenomena arising when
comparing Woronowicz algebras with Popa systems, and will
be quoted along the rest of the paper.

\lvs
\noindent {\bf Example 1.1} If $A$ is a Woronowicz algebra then the
following conditions are equivalent (cf. the formulas {\it
(f1-4)} see \cite{w1}, \cite{bs}):

- the Haar measure $h\in A^*$ is a trace.

- the square of the antipode $\kappa$ is the identity of $A_s$.

- all $f_z$'s are equal to the counit.

This happens for instance if $A$ is of the form $C(G)$ with $G$
compact group, or $\c^*(\Gamma )$ with $\Gamma$ discrete
group, or when $A$ is finite dimensional. In this case it is
clear that
$Q_v=id_H$, $d(v)=dim(H)$, and $\hat{v}=\vb$ for any unitary
corepresentation $v\in\l (H)\otimes A_s$.

\noindent {\bf Example 1.2} The $q$-deformation of the $\c^*$-algebras
$C(G)$, with $G$ compact classical Lie group and $q>0$ was
started by Woronowicz, and completed by Rosso \cite{ro}, via a
beautiful application of the Tannakian duality \cite{w2}. By
\cite{r1},\cite{ro} the fusion semiring is invariant
under $q$-deformations.

If $\mu\in [-1,1]-\{ 0\}$ and $u$ is the fundamental
representation of $\s_\mu\u (2)$ then
$Q_u=diag(\mid\mu\mid^{-1},\mid\mu\mid )$ (cf. the Appendix of
\cite{w1}).

\noindent {\bf Example 1.3} For $n\in\n$ and $F\in\g\ll (n,\c )$
satisfying $F\overline{F}\in\r I_n$ let
$A_o(F)$ be the universal $\c^*$-algebra generated by the
entries of a $n\times n$ matrix $u$ with the relations
$u=F\overline{u}F^{-1}=$ unitary. This Woronowicz algebra was
introduced in \cite{vdw} as to represent the ``free analogue
of $\o (n)$'' and its fusion semiring was shown in \cite{b0}
to be isomorphic to the one of $\s\u (2)$. Moreover, by
combining this result with \cite{kw} one finds that the
category of corepresentations of $A_o(F)$ is always monoidal
equivalent to the category of corepresentations of $\s_\mu\u
(2)$ for some $\mu\in [-1,1]-\{ 0\}$. This is a relevant
example for the considerations in the end of the introduction. 

We mention that one can show that the Woronowicz algebras
$A_o(F)$ are exactly the ``compact forms'' of Gurevich's
quantizations of $\s\ll (2,\c )$ \cite{g}.

{\noindent \bf Example 1.4} For $n\in\n$ and $F\in\g\ll (n,\c )$ let
$A_u(F)$ be the universal $\c^*$-algebra generated by the
entries of a $n\times n$ matrix $u$ with the relations $u$
unitary, $F\overline{u}F^{-1}$ unitary. These
Woronowicz algebras correspond to both notions of
``algebras of continuous functions on the free unitary groups''
and ``$\c^*$-algebras of the free free groups''
(sic!) \cite{vdw}. The fusion semiring of $A_u(F)$ was computed
in \cite{b1}. Using {\it (ii) $\Rightarrow$ (i)} in lemma
1.2 it's easy to see that $Q_u$ is a scalar multiple of the transpose of 
$\sqrt{F^*F}$.

Thus if $n\in\n$ and $Q\in M_n(\c )$ is a positive diagonal
matrix such that $Tr(Q^2)=Tr(Q^{-2})$, then for $(A_u(Q),u)$
we have $Q_u=Q$ (note: one can show that given any $n\in\n$ and
$F\in\g\ll (n,\c )$ there exists such a matrix $Q$ such
that $A_u(F)$ is isomorphic to $A_u(Q)$, and that $Q$ is
unique up to $Q\mapsto Q^{-1}$ and up to the permutation of
its diagonal entries; this classification of the algebras
$A_u(F)$ could be easily deduced from \cite{b1} and from the above computation of $Q_u$).

\begin{lemm}
(i) $Q_{\hat{v}}=(Q_v^t)^{-1}$, $\hat{\hat{v}} =v$ and $d(\hat{v})=d(v)$.

(ii) $Q_{v\otimes w}=Q_v\otimes Q_w$, $\widehat{v\otimes w}=\hat{w}\otimes\hat{v}$ and $d(v\otimes w)=d(v)d(w)$.

(iii) $Q_{v+w}=diag(Q_v,Q_w)$, $\widehat{v+w}=\hat{v}+\hat{w}$ and $d(v+w)=d(v)+d(w)$.
\end{lemm}

{\em Proof.} The formulas for the quantum dimensions and for the
canonical duals follow easily from the formulas for the $Q$'s,
which we will now prove. First,
$Q_{v+w}=diag(Q_v,Q_w)$ is clear and $Q_{v\otimes w}=Q_v\otimes
Q_w$ follows from the fact that $f_{\frac{1}{2}}$ is a
character. Also using {\it (f1)} and {\it (f2)} we get that 

\lvs\noindent $Q_{\hat{v}}=(id\otimes f_{\frac{1}{2}})(t\otimes
f_{\frac{1}{2}}*\kappa (.)*f_{-\frac{1}{2}})v=(t\otimes
f_{\frac{1}{2}}(f_{\frac{1}{2}}*\kappa
(.)*f_{-\frac{1}{2}}))v=$

\lvs\noindent $(t\otimes
(f_{\frac{1}{2}}*f_{\frac{1}{2}}*f_{-\frac{1}{2}})\kappa )v=
(t\otimes f_{\frac{1}{2}}\kappa
)v=(t\otimes f_{-\frac{1}{2}})v=(Q_v^t)^{-1}\,\,\,\qed$ 

\lvs
If $v$ is a finite dimensional unitary corepresentation of
$A_s$ then the contragradient $v^c$, the complex conjugate
$\vb$ and the canonical dual $\hat{v}$ are of course in the same
equivalence class, which is dual to the equivalence class of
$v$. In the next proposition we
state the precise (spatial) form of duality between $v$ and
$\hat{v}$ (see the remark below).
\lvs
\noindent {\bf Notation.} If $H$ is a finite dimensional Hilbert
space $\xi :\l (H)\simeq H\otimes\Hb$ is the canonical linear
map.

\begin{prop} 
For any finite dimensional unitary corepresentation $v\in\l
(H)\otimes A_s$ define $i_v\in \l (\c ,H\otimes\Hb )$ by
$1\mapsto\xi (Q_v)$. Define also $p_v:=i_{\hat{v}}^*$. Then the
following  ``duality formulas'' hold:

{\it (d1)} $i_v\in 
Hom(1,v\otimes \hat{v})$ and $p_v\in Hom(\hat{v}\otimes v,1)$.

{\it (d2)} $(id_v\otimes p_v)(i_v\otimes id_v)=id_v$ and
$(p_v\otimes id_{\hat{v}})(id_{\hat{v}}\otimes
i_v)=id_{\hat{v}}$.

{\it (d3)} $p_{\hat{v}}=i_v^*$, $p_v=i_{\hat{v}}^*$
and $p_vi_{\hat{v}}=p_{\hat{v}}i_v=d(v)=d(\hat{v})$.

{\it (d4)} $p_{r\otimes s}=p_s(id_{\hat{s}}\otimes p_r\otimes
id_s)$ and 
$i_{r\otimes s}=(id_r\otimes i_s\otimes id_{\hat{r}})i_r$,
$\forall\, r,s$.
\end{prop}

{\em Proof.} {\it (d3)} is trivial, and 
{\it (d4)} follows from lemma 1.3 (ii). Let us choose an
orthogonal basis $\{ e_1,...,e_n\}$ of $H$ consisting of 
eigenvectors of $Q_v$, so that $Q_v=\sum Q_{ii}e_{ii}$.
Then $<(id_v\otimes p_v)(i_v\otimes id_v)x,y>$ is equal to 
$$<\xi
(Q_v)\otimes x,y\otimes \xi (Q_{\hat{v}})>=\sum
<Q_{ii}e_i\otimes
\overline{e}_i\otimes x,y\otimes Q_{jj}^{-1}\otimes
\overline{e}_j\otimes e_j>=<x,y>$$
for any $x,y\in H$, and {\it (d2)} follows. It remains to
verify that $i_v\in  Hom(1,v\otimes \hat{v})$ for any $v$; the
other assertion $p_v\in Hom(\hat{v}\otimes v,1)$ will follow
from $p_v=i_{\hat{v}}^*$ and from $i_{\hat{v}}\in
Hom(1,\hat{v}\otimes v)$. We have $i_v(1)=\sum Q_{ii}e_i\otimes
\overline{e}_i$, and if $\zeta :=\sum e_i\otimes
\overline{e}_i$ then
$$(v\otimes \vb )(\zeta\otimes 1)=\sum_{i,k,a} e_i\otimes \overline{e}_k\otimes v_{ia}v_{ka}^*=\sum_{i,k} 
e_i\otimes \overline{e}_k\otimes \delta_{ik}1=(\zeta\otimes 1)$$
As $\vb =(Q_v^{-1})^t\hat{v}(Q_v)^t$ we get 
$(v\otimes \hat{v})(\sum Q_{ii} e_i\otimes \overline{e}_i\otimes
1)=(\sum Q_{ii} e_i\otimes \overline{e}_i\otimes 1)$,  i.e.
$i_v\in Hom(1,v\otimes \hat{v})$. $\qed$

\lvs
\noindent {\bf Remark.} Recall that in a monoidal category a  duality
$(X,Y,i_X,p_X)$  consists of two objects $X,Y$ and two
morphisms $i_X\in Hom(1,X\otimes Y)$, $p_X\in Hom(Y\otimes
X,1)$  such that the compositions
$$\displaystyle{X=1\otimes
X\mathop{\longrightarrow}^{i_X\otimes id} X\otimes Y\otimes X 
\mathop{\longrightarrow}^{id\otimes p_X}X\otimes 1=X}$$ 
$$\displaystyle{Y=Y\otimes 1\mathop
{\longrightarrow}^{id\otimes i_X} Y\otimes X\otimes Y 
\mathop{\longrightarrow}^{p_X\otimes id}1\otimes Y=Y}$$
are equal to $id_X$ and $id_Y$ correspondingly (in this case
Frobenius reciprocity holds, etc. - see for instance \cite{kw}).
Thus {\it (d1,d2)} show that
$(v,\hat{v},i_v,p_v)$ and
$(\hat{v},v,i_{\hat{v}},p_{\hat{v}})$ are dualities, i.e. that
$\hat{v}$ is both a right and a left dual for $v$ is the
monoidal category of finite dimensional corepresentations of
$A$.
\section{Algebras of symmetries}

We use the following notation for tensor products and related lattices
(we recall that in this paper the word {\em lattice} will always mean
{\em system of inclusions}; see the introduction). 

\noindent {\bf Notations.} Let $\n^{*2}=<\alpha ,\beta >$ be the free
monoid consisting of words in $\alpha ,\beta$, plus the empty
word denoted by $e$. We define an involution
$\,\widehat{ }\,$ on $\n^{*2}$ by $\hat{\alpha}=\beta$,
$\hat{\beta}=\alpha$ and antimultiplicativity. 

If $({\cal C},\otimes ,1,\,\widehat{ }\,)$ is an involutive
monoid, and if $A\in {\cal C}$, $x\in \n^{*2}$ we
denote by $A^{\otimes x}$ the image of $x$ by the unique
morphism $\n^{*2}\rightarrow {\cal C}$ sending
$\alpha\mapsto A$. 

If $0\leq i\leq j$ 
are integers we define $[i,j]\in\n^{*2}$ to be
$(\alpha\beta )^{\frac{j-i}{2}}$ if $i,j$ are even,
$(\beta\alpha )^{\frac{j-i}{2}}$ if $i,j$ are odd,
$(\alpha\beta )^{\frac{j-i-1}{2}}\alpha$ if $i$ is even and $j$
is odd, and $\beta (\alpha\beta )^{\frac{j-i-1}{2}}$ if $i$ is
odd and $j$ is even.

Examples: $[2,6]=\alpha\beta\alpha\beta$, $[1,4]=\beta\alpha\beta$,
$[3,3]=e$, so that $A^{[2,6]}=A\otimes\hat{A}\otimes A\otimes\hat{A}$,
$B^{[1,4]}=\hat{B}\otimes B\otimes\hat{B}$ and $C^{[3,3]}=1$ if
$A,B,C$ are elements of an involutive monoid $({\cal C},\otimes ,1,\,\widehat{ }\,)$.

\lvs
We will use these notations when $({\cal
C},\otimes ,1)$ is the monoidal category of unitary
corepresentations of a Woronowicz algebra, with $\,\widehat{
}\,=$ the canonical dual, or when
$({\cal C},\otimes ,1)$ is the monoidal category of complex
vector spaces, with $\,\widehat{ }\,=$ complex conjugation.
Thus the lattices $L(\pi )$ and $L(H)$ in the Introduction are
equal to $End(\pi^{\otimes [i,j]})_{0\leq
i\leq j <\infty}$ and $End(H^{\otimes [i,j]})_{0\leq
i\leq j <\infty}$ respectively.

\lvs 
Recall that a standard $\lambda$-lattice of commuting squares
(or {\em Popa system}) is a system $(A_{ij})_{0\leq i\leq j
<\infty}$ of finite dimensional 
$\c^*$-algebras with $A_{ii}=\c$, $A_{ij}\subset A_{kl}$ for $k\leq i,j\leq l$, and with a given faithful trace $\tau$ 
on $\cup_{i,j}A_{ij}$ satisfying the following properties 
(see section 1 of \cite{p3}):

{\em The commuting square condition.} $E_{A_{ij}}E_{A_{kl}}=E_{A_{kl}}E_{A_{ij}}=E_{A_{rs}}$ 
where $r=max\{ i,k\}$, $s=min\{ j,l\}$ and $E_B$ is the $\tau$-preserving conditional expectation onto $B$.

{\em Jones conditions.} There exists a
representation of the $\lambda$-sequence of Jones projections
$\{ e_i\}_{i\geq 2}$ in $\cup_{i,j} A_{ij}$ such that $e_j\in
A_{i-2,k}$, $\forall\, 2\leq i\leq j\leq k$ and

$e_{j+1}xe_{j+1}=E_{A_{i,j-1}}(x)e_{j+1}$, 
$\forall\, i\leq j-1$ and $x\in A_{ij}$.

$e_{i+1}xe_{i+1}=E_{A_{i+1,j}}(x)e_{i+1}$, 
$\forall\, i+1\leq j$ and $x\in A_{ij}$.

{\em Markov conditions.} 

$\lambda^{-1}E_{A_{i,j+1}}(xe_{j+2})e_{j+2}=xe_{j+2}$, 
$\forall\, j\geq i$ and $x\in A_{i,j+2}$.

$\lambda^{-1}E_{A_{i+1,j}}(xe_{i+2})e_{i+2}=xe_{i+2}$, 
$\forall\, j\geq i+2$ and $x\in A_{i,j}$.

{\em Commutation relations.} $[A_{ij},A_{kl}]=0$, $\forall\, i\leq j\leq k\leq
l$.

More precisely, by {\em Popa system} we mean the triple consisting of the
lattice, the trace, and the sequence of Jones projections, which
satisfies the above four axioms. Two
Popa systems $((A_{ij})_{0\leq i\leq j<\infty},\tau
,(e_k)_{k\geq 2})$ and $((B_{ij} )_{0\leq i\leq
j<\infty},\varphi ,(f_k )_{k\geq 2})$ are said to
be equal if there exists a system of trace-preserving,
Jones projections-preserving and inclusion-preserving
$\c^*$-isomorphisms $A_{ij}\rightarrow B_{ij}$. 

\begin{theo}
If $v$ is a finite dimensional unitary corepresentation of
a Woronowicz algebra then the lattice $End(v^{\otimes
[i,j]})_{0\leq i\leq j <\infty}$ (together with the canonical
traces from the definition 1.1 and with the Jones projections given by
$$e_{2k}:=d(v)^{-1}i_vp_{\hat{v}}\in
End(v\otimes\hat{v}),\,\,\,\,\,
e_{2k+1}:=d(v)^{-1}i_{\hat{v}}p_v\in End(\hat{v}\otimes v).$$
for any $k\geq 1$) is a Popa $\lambda$-system, with
$\lambda =d(v)^{-2}$.
\end{theo}

{\em Proof.} We have to verify the above four axioms. The
Commutation relations are clear. Note that all the inclusions
between the algebras in the lattice are of the form
$$End(a)\hookrightarrow End(r\otimes a\otimes w),\,\, x\mapsto id_r\otimes x\otimes id_w$$
Let us define a linear map 
$E_{r,a,w}:End(r\otimes a\otimes w)\rightarrow End(a)$ by 
$$x\mapsto d(r)^{-1}d(w)^{-1}
(p_r\otimes id_a\otimes p_{\hat{w}})(id_{\hat{r}}\otimes
x\otimes id_{\hat{w}}) (i_{\hat{r}}\otimes id_a\otimes i_w)$$
By using {\it (d4)} we get 
$$E_{r,a,w}E_{r^\prime ,r\otimes a\otimes w,w^\prime}=E_{r^\prime\otimes r,a,w\otimes w^\prime}\,\,\,\, \forall\, r^\prime ,r,a,w
,w^\prime\,\,\,\,\, (\star )$$ 
We verify now that $E_{r,a,w}$ is the conditional expectation, i.e. that 
$$E_{r,a,w}(id_r\otimes x\otimes id_w)=x$$
$$E_{r,a,w}[(id_r\otimes x\otimes id_w)y(id_r\otimes z\otimes id_w)]=xE_{r,a,w}(y)z$$
$$\tau_{r\otimes a\otimes w}=\tau_a E_{r,a,w}$$
$\forall\, x,z\in End(a)$ and $y\in End(r\otimes a\otimes w)$. 
Indeed, the first formula follows from {\it (d3)}, the second
one is clear by definition of $E_{r,a,w}$, and the third one
follows from $\tau_v=E_{1,1,v}=E_{v,1,1}$ for any $v$ and from 
$(\star )$. The commuting square axiom follows now from 
$(\star )$. 

By {\it (d3)}, the $e_s$'s are selfadjoint projections. Using {\it (d2)} it's easy to see that $\{ e_2,e_3,...\}$ is a
representation of the $\lambda$-sequence of Jones projections, 
with $\lambda :=d(v)^{-2}$. We will prove the other two axioms
for the  horizontal inclusions (same proof for the vertical
ones). Note that every pair of consecutive horizontal
inclusions is of the form 
$$End(a)\hookrightarrow End(a\otimes w)\hookrightarrow End(a\otimes w\otimes\hat{w})$$
with $w=v$ or $\hat{v}$, and the Jones projection in 
$End(a\otimes w\otimes\hat{w})$ is $e:=d(w)^{-1}id_a\otimes
i_wp_{\hat{w}}$. The Jones condition, namely
$$e(x\otimes id_{\hat{w}})e=(E_{1,a,w}(x)\otimes id_{w\otimes\hat{w}})e,\,\, \forall\, x\in End(a\otimes w)$$
is trivial from the definitions of $e$ and $E_{1,a,w}$. 
The Markov condition follows by applying {\it (d2,d3)} several
times  - if $x\in End(a\otimes w\otimes\hat{w})$ then 
\lvs 
$d(w)^2(E_{1,a\otimes w,\hat{w}}(xe)\otimes id_{\hat{w}})e=$

$[(id_{a\otimes w}\otimes p_w)(x\otimes id_w)(id_a\otimes
i_w\otimes id_w)\otimes id_{\hat{w}}]e=$

$d(w)^{-1}(id_{a\otimes w}\otimes p_w\otimes
id_{\hat{w}})(x\otimes id_{w\otimes \hat{w}})(id_a\otimes
i_w\otimes i_wp_{\hat{w}})=$

$d(w)^{-1}(id_{a\otimes w}\otimes [(p_w\otimes
id_{\hat{w}})(id_{\hat{w}}\otimes i_wp_{\hat{w}})])
([x(id_a\otimes i_w)]\otimes id_{w\otimes\hat{w}})=$

$d(w)^{-1}(id_{a\otimes w\otimes\hat{w}}\otimes
p_{\hat{w}})([x(id_a\otimes i_w)]\otimes id_{w\otimes\hat{w}})=$

$d(w)^{-1}(x\otimes id_{w\otimes \hat{w}})(id_a\otimes
i_w\otimes p_{\hat{w}})=xe$. $\qed$
\lvs

The Popa systems are exactly the lattices of $\c^*$-algebras
(with traces) which can arise as lattices of higher relative
commutants of extremal inclusions of finite index $II_1$
factors \cite{p3}. There are at least two kind of examples when
the Popa system $End(v^{\otimes [i,j]})_{0\leq i\leq j
<\infty}$ is the one corresponding to some ``nice'' subfactor
associated to
$(A,v)$:
\lvs 
\noindent {\bf Example 2.1 (cf. \cite{was})} Let $G\rightarrow Aut(P)$
be a minimal action of a compact group $G$ on a $II_1$ factor
$P$. If $\pi :G\rightarrow \l (H)$ is a finite dimensional
unitary representation of $G$ then the lattice of higher
relative commutants of the inclusion of fixed point algebras
$$P^G\subset (P\otimes\l (H))^G$$
is $End(v^{\otimes [i,j]})_{0\leq i\leq j <\infty}$,
where $v$ is the corepresentation corresponding to
$\pi$ of the Woronowicz algebra $A=C(G)$.

\lvs 
\noindent {\bf Example 2.2 (cf. \cite{p1})} Let $\Gamma\subset Aut(P)$ be
an outer (discrete) group of automorphisms of a $II_1$ factor
$P$. If $\Gamma$ is generated by $g_1,...,g_n$ then 
the lattice of higher relative commutants of the inclusion 
$$P\subset M_n(P)$$
given by $x\mapsto diag(g_1(x),...,g_n(x))$ is $End(v^{\otimes
[i,j]})_{0\leq i\leq j <\infty}$, where $v$ is the
corepresentation $diag(u_{g_1},...,u_{g_n})$ of the Woronowicz
algebra $A=\c^*(\Gamma )$.

\lvs
It is possible to extend the above results to coactions/actions
of more general Woronowicz algebras (see for instance
\cite{pw}), and in fact this quantum group formalism shows
that the above two constructions are of the same nature (i.e.
one can pass from each of them to the other one by considering
dual coactions/actions; I owe this observation from E.
Blanchard). In the general case Popa's universal construction
is available:
\lvs 
\noindent {\bf Example 2.3} Let $\mu\in (0,1)$ and consider the
fundamental representation $v$ of $\s_\mu\u (2)$. As
$d(v)=\mu^2+\mu^{-2}$ and $End(v)\simeq\c$ (see example
1.2), any inclusion having $End(v^{\otimes [i,j]})_{0\leq i\leq
j <\infty}$ as lattice of higher relative commutants has to be
irreducible and of index $\lambda^{-1}=(\mu^2+\mu^{-2})^2>4$.

The algebras $End(v^{\otimes [i,j]})$ being generated by the
Jones projections (cf. the representation theory of
$\s_\mu\u (2)$, see \cite{w2}), $End(v^{\otimes
[i,j]})_{0\leq i\leq j <\infty}$ is the lattice of higher
relative commutants of the subfactors constructed in \cite{p2}.

\section{Representations of Popa systems}

There are several restrictions on the Popa systems associated
to corepresentations of Woronowicz algebras, for instance the
index is always $\geq 4$ and the index of the amenable ones is
the square of an integer (see section 6). We will prove that
in fact the only extra structure of a Popa system of the form
$End(v^{\otimes [i,j]})_{0\leq i\leq j <\infty}$ comes from the
fact that it can be embedded into the lattice $\l (H^{\otimes
[i,j]})_{0\leq i\leq j<\infty}$, where $H$ is the Hilbert space
where the corepresentation $v$ acts.   
\lvs
\noindent {\bf Definition 3.1} A representation of a Popa system $(A_{ij})_{0\leq i\leq
j<\infty}$ on a finite dimensional Hilbert space $H$ is a
system $\pi =(\pi_{ij})_{0\leq i\leq
j<\infty}$ of inclusions of $\c^*$-algebras
$\pi_{ij}:A_{ij}\hookrightarrow
\l (H^{\otimes [i,j]})$ such that
$$\pi_{kl}(x)=id_{H^{\otimes [k,i]}}\otimes\pi_{ij}(x)\otimes
id_{H^{\otimes [j,l]}}$$
for every $0\leq k\leq i\leq j\leq l$ and $x\in A_{ij}$. 
\lvs

\noindent {\bf Example 3.1} If $A$ is a Woronowicz algebra and $v
\in\l (H)\otimes A_s$ is a finite dimensional unitary
corepresentation of $A$ then $v^{\otimes x}$ acts on
$H^{\otimes x}$ for every $x\in \n^{*2}$, so $End(v^{\otimes
x})$ is by definition a subalgebra of $\l (H^{\otimes x})$. It
follows that
$End(v^{\otimes [i,j]})_{0\leq i\leq j <\infty}$ is canonically
represented on
$H$.
\lvs 
\noindent {\bf Example 3.2} Fix $\lambda \in (0,1/4]$ and let
$(A_{ij}^\lambda )_{0\leq i\leq j<\infty}$ be the Popa system
generated by the $\lambda$-sequence of Jones projections (i.e.
the one in example 2.3). If
$H$ is a finite dimensional Hilbert space and $Q\in\l (H)$ is
a positive operator satisfying
$Tr(Q^2)=Tr(Q^{-2})=\lambda^{-1/2}$, it is easy to see that
there exists a (unique) representation of $(A_{ij}^\lambda
)_{0\leq i\leq j<\infty}$ on $H$, sending 
$\pi_{2k-2,2k}: e_{2k}\mapsto Proj_{\c\,\xi (Q)}$
and $\pi_{2k-1,2k+1}: e_{2k+1}\mapsto Proj_{\c\,\xi
((Q^{-1})^t)}$. 

In fact these are (modulo some ``unitary
equivalence'', see below) all the representations of
$(A_{ij}^\lambda )_{0\leq i\leq j<\infty}$ on finite
dimensional Hilbert spaces. Note that the so-called PPTL representation (see
\cite{pipo} and section 5 in \cite{jcbms}) is a particular case of this construction.
\lvs 
The second example is in fact a particular case of the first
one, for $v=$ the fundamental representation of $A_u(Q)$. In
fact the following result holds:

\begin{theo}
For any Popa system $(A_{ij})_{0\leq i\leq j<\infty}$ the
following conditions are equivalent:

(i) There exists a finite dimensional unitary
corepresentation $v$ of a Woronowicz algebra $A$ such that
$(A_{ij})_{0\leq i\leq j<\infty}=End(v^{\otimes
[i,j]})_{0\leq i\leq j <\infty}$.

(ii) There exists a finite dimensional Hilbert space $H$ and a
representation $\pi$ of $(A_{ij})_{0\leq i\leq j<\infty}$ on
$H$ such that $\pi_{02}(e_2)$ is a
rank one projection.
\end{theo}

The implication {\it (i)} $\Longrightarrow$ {\it (ii)} is
clear, and the converse will be proved in the next
section. In the rest of this section we give some preliminary
results on the representations of Popa systems.

Recall that the Jones projections of $End(v^{\otimes
[i,j]})_{0\leq i\leq j <\infty}$ have a special form - they
are the rank one projections onto $\c\,\xi (Q_v)$ or onto
$\c\,\xi ((Q_v^{-1})^t)$ (where $\xi :\l (H)\simeq
H\otimes\Hb$ is the canonical linear map). We will need the
following proposition on ``unitary
equivalence'' of representations.
\lvs
\noindent {\bf Definition 3.2} A representation $\pi$ of a Popa system $(A_{ij})_{0\leq i\leq
j<\infty}$ on $H$ is said to be normalised if there exists a
positive operator $Q\in\l (H)$ satisfying
$Tr(Q^2)=Tr(Q^{-2})=\lambda^{-1/2}$ such that $\pi_{2k-2,2k}
(e_{2k})$ is the projection onto
$\c\,\xi (Q)$ and $\pi_{2k-1,2k+1} (e_{2k+1})$ is the
projection onto $\c\,\xi ((Q^{-1})^t)$, $\forall\, k\geq 1$.

\begin{prop}
If $\pi$ is a representation such that $\pi_{02}(e_2)$ is a
rank one projection, then there exist two sequences of unitaries
$U_1,U_3,U_5,...\in \l (H)$ and $U_2,U_4,U_6,...\in\l (\Hb )$ such that 
$$\pi^\prime_{ij}=
ad(U_{i+1}\otimes U_{i+2}\otimes ...\otimes U_j)\circ 
\pi_{ij}$$
is normalised representation.
\end{prop}

We will use the following elementary lemma:

\begin{lemm}
Let $H$ be a finite dimensional Hilbert space and $e\in \l
(H\otimes\Hb )$, $f\in \l (\Hb\otimes H)$ be rank one
projections satisfying the Jones relation
$$(id_H\otimes f)(e\otimes id_H)
(id_H\otimes f)=\lambda (id_H\otimes f)$$
for some $\lambda >0$. If $E,F\in\l (H)$ are such
that $\xi (E)$ and $\xi (F^t)$ are norm one vectors in the
images of $e$ and $f$, then $(F^*F)(E^*E)=\lambda id_H$. 
\end{lemm}

{\it Proof.} If $\{ h_k\}$ is an orthogonal
basis of $H$ then a straighforward computation
shows that
$$<(id\otimes f)(e\otimes id)(id\otimes f)(h_a\otimes \xi (F^t)
),(h_b\otimes\xi (F^t) )>=(EF^*FE^*)_{ba}$$ 
for any $a,b$, so the Jones relation
implies $EF^*FE^*=\lambda I_H$. $\qed$ 

\lvs 
{\it Proof of proposition 3.1.} Note first that the Jones
relations imply (by induction on $n\geq 2$) that each
$\pi_{n-2,n}(e_n)$ is a rank one projection. 
\lvs 
{\it Step I. Construction of $Q,U_1,U_2$.} The above
lemma applies with $e=\pi_{02}(e_2)$ and $f=\pi_{13}(e_3)$. 
Thus if $e$ is the projection onto $\c\,\xi (E)$ then $E$ is 
invertible; let $E=QU$ be its polar decomposition. If 
$\{ h_1,...,h_n\}$ is an orthogonal basis of
$H$ consisting of eigenvectors of $Q$, then $e$ is the
orthogonal projection onto $\c\,\sum
Q_{ii}U_{ij}h_i\otimes\overline{h_j}$. By multiplying $Q$ with
a scalar, we may suppose that
$Tr(Q^2)=Tr(Q^{-2})$. We have 
$$(I_n\otimes \overline{U})(\sum Q_{ii}U_{ij}h_i\otimes\overline{h_j})=\sum Q_{ii}h_i\otimes \overline{h}_i =\xi (Q)$$
so $(I_n\otimes \overline{U})e(I_n\otimes U^t)$ is the
orthogonal projection onto $\c\,\xi (Q)$. By perturbing the
representation $\pi$ with the sequence of unitaries $\{ I_n,
\overline{U},I_n,I_n,...\}$ we may assume that
$e=\pi_{02}(e_2)$ is the orthogonal projection onto $\c\,\xi
(Q)$.
\lvs
{\it Step II. Construction of $U_s$, $s\geq 3$ and proof of
$Tr(Q^2)=\lambda^{-1/2}$.} Define
$e_{2k}^\prime$ to be the projection onto $\c\,\xi (Q)$ and
$e_{2k+1}^\prime$ to be the projection onto $\c\,\xi (Q^{-1})$,
for every $k\geq 1$. We have to show that if
$\pi_{k-2,k}(e_k)=e_k^\prime$ for
$k=2,3,...,s-1$ then $(id\otimes
ad(U_s))(\pi_{s-2,s}(e_s))=e_s^\prime$ for some
unitary $U_s$. It is enough to do it for $s=3$. Let
$e=\pi_{02}(e_2)$ and
$f=\pi_{13}(e_3)$ and assume that $f$ is the orthogonal
projection onto $\c\,\nu$, with $\nu =\xi (F^t)$ a norm one
vector. The above lemma shows that 
$$(F^*F)(Tr(Q^2)^{-1}Q^2)=\lambda I_n$$
and as $\nu$ is a norm one vector, $Tr(F^*F)=1$,
so $Tr(Q^{-2})=Tr(Q^2)=\lambda^{-1/2}$. Also
$F^*F=\lambda^{1/2}Q^{-2}$, so $F=\lambda^{1/4}UQ^{-1}$ for
some unitary $U$. Thus $f$ is the orthogonal projection onto
$\c\,\sum Q_{ii}^{-1}U_{ji}\overline{h}_i\otimes h_j$. As
$(I_n\otimes U^*)\sum Q_{ii}^{-1}U_{ji}\overline{h}_i\otimes
h_j =\sum Q_{ii}^{-1}\overline{h}_i\otimes h_i$, one can choose
$U_3:=U$. $\qed$

\lvs
Recall that the conditional expectations in a Popa system of
the form $End(v^{\otimes [i,j]})_{0\leq i\leq j <\infty}$ have
a certain special form; we will prove in the next proposition
that the same happens for a Popa system having a normalised
representation. 

\begin{prop}
Let $(A_{ij})_{0\leq i\leq j<\infty}$, $H$, $\pi$, $Q$ be as
in definition 3.2. Define linear maps
$i_{\alpha}\in \l (\c ,H\otimes\Hb )$ by $1\mapsto\xi (Q)$
and $i_{\beta}\in \l (\c, \Hb\otimes H)$ by
$1\mapsto\xi ((Q^{-1})^t)$. Let
$p_{\alpha}:=i_{\beta}^*$ and
$p_{\beta}:=i_{\alpha}^*$. For
$x\in\n^{*2}$ let $id_x$ be the identity of $H^{\otimes x}$. 

(i) The following ``duality formulas''
hold, $\forall\,\gamma\in\{\alpha ,\beta\}$:
$$(id_\gamma\otimes p_\gamma )(i_\gamma\otimes id_\gamma
)=(p_{\hat{\gamma}}\otimes id_\gamma )(id_\gamma \otimes
i_{\hat{\gamma}})=id_\gamma$$

(ii) The Jones projections of $(A_{ij})_{0\leq i\leq
j<\infty}$ are given by 
$$\pi_{2k-2,2k}(e_{2k})=\lambda^{1/2}i_{\alpha}p_{\beta}
,\,\,\,
\pi_{2k-1,2k+1} (e_{2k+1})=\lambda^{1/2}i_{\beta}
p_{\alpha}$$

(iii) The horizontal conditional expectations are given by
$$\pi_{i,j-1}(E_{A_{i,j-1}}(T))=\lambda^{1/2}
(id_{[i,j-1]}\otimes
p_{\hat{\gamma}} ) (\pi_{ij}(T)\otimes id_{\hat{\gamma}}
)(id_{[i,j-1]}\otimes i_\gamma )$$
$\forall\, 0\leq i\leq j$ and $T\in A_{ij}$, where $\gamma
=\alpha$ if $j$ is odd and $\gamma
=\beta$ if $j$ is even. A similar formula holds for the
vertical conditional expectations.
\end{prop}

{\it Proof.} (ii) is trivial, for (i) see the
proof of proposition 1.1 {\it (d2)}. For
(iii) recall that by the {\em Jones conditions} the
following equality holds in $A_{i,j+1}$:
$$E_{A_{i,j-1}}(T)e_{j+1}=e_{j+1}Te_{j+1}$$
As $\pi_{i,j+1}(e_{j+1})=\lambda^{1/2}id_{[i,j-1]}\otimes
i_\gamma p_{\hat{\gamma}}$, by applying $\pi_{i,j+1}$ the Jones
condition becomes
$$\pi_{i,j-1}(E_{A_{i,j-1}}(T))\otimes (i_\gamma
p_{\hat{\gamma}} )=
\lambda^{1/2}(id_{[i,j-1]}\otimes i_\gamma
p_{\hat{\gamma}} )(\pi_{ij}(T)\otimes
id_{\hat{\gamma}})(id_{[i,j-1]}\otimes i_\gamma
p_{\hat{\gamma}} )$$ 
which is equal to $\lambda^{1/2}
(id_{[i,j-1]}\otimes
p_{\hat{\gamma}} ) (\pi_{ij}(T)\otimes id_{\hat{\gamma}}
)(id_{[i,j-1]}\otimes i_\gamma )\otimes (i_\gamma
p_{\hat{\gamma}} )$. $\qed$

\lvs
\noindent {\bf Remark.} Note that only the Jones conditions were used
in the above proof. In fact one can prove that any lattice
$(A_{ij})_{0\leq i\leq j<\infty}$ of $\c^*$-algebras (with
traces) which satisfies the {\it Jones conditions} and which
has a normalised representation is a Popa system (the
{\em Commutation relations} are trivial, the {\em Commuting
square condition} is clear from the above formula for the
conditional expectations, and the {\em Markov conditions} may
be checked by the same computation as in the proof of theorem 2.1). This is the reason why in the next
section we will use proposition 3.2 instead of Popa's four
axioms. Actually, we will need not only proposition 3.2,
but also its corollaries 3.1 and 3.2 below.

\begin{coro}
If $0\leq i\leq j$ and $T\in A_{ij}$ then 
$$\pi_{i,j-2}(E_{A_{i,j-2}}(e_jTe_j))=(id_{[i,j-2]}\otimes
p_{\hat{\gamma}} )\pi_{ij}(T)(id_{[i,j-2]}\otimes i_\gamma )$$
where $\gamma =\beta$ if $j$ is odd and $\gamma 
=\alpha$ if $j$ is even. 
\end{coro}

{\it Proof.} Follows from
$E_{A_{i,j-2}}=E_{A_{i,j-2}}E_{A_{i,j-1}}$ and by applying
twice proposition 3.2. (iii). $\qed$

\begin{coro}
If $0\leq i\leq j$ and $0\leq k\leq l$ are such that
$[i,j]=[k,l]$ in $\n^{*2}$, then
$\pi_{i,j}(A_{i,j})=\pi_{k,l}(A_{k,l})$.
\end{coro}

{\it Proof.} By an induction argument, it is enough to prove it
for $k=i+2$ and $l=j+2$. Recall that the canonical
isomorphism $sh:A_{ij}\rightarrow A_{i+2,j+2}$ satisfies the
formula (see for instance \cite{bi})
$$sh(T)e_{i+2}=\lambda^{i-j}e_{i+2}e_{i+3}...e_{j+1}
Te_{j+2}e_{j+1}...e_{i+3}e_{i+2}\,\,\,\,\,\, (\star )$$
Now by proposition 3.2 (i,ii) we get that 
$$\pi_{i,j+2}(e_{i+2}...e_{j+1})=
\lambda^{j-i/2}i_{\hat{\gamma}}\otimes id_{[i,j-1]}\otimes
p_\delta\otimes id_{\hat{\delta}}$$
$$\pi_{i,j+2}(e_{j+2}...e_{i+2})=\lambda^{j-i+1/2}p_\gamma
\otimes id_{[i,j]}\otimes i_\delta$$
where $\gamma ,\delta\in\{ \alpha ,\beta\}$ depend on the
parity of $i$ and $j$. Thus by applying $\pi_{i,j+2}$ to
$(\star )$ we get
that $\lambda^{1/2}i_{\hat{\gamma}}p_\gamma\otimes
\pi_{i+2,j+2}(sh(T))$
is equal to 
$$\lambda^{1/2}i_{\hat{\gamma}}p_\gamma\otimes 
(id_{[i,j-1]}\otimes p_\delta\otimes
id_{\hat{\delta}})(\pi_{ij}(T)\otimes
id_{\delta\hat{\delta}})(id_{[i,j]}\otimes i_\delta )$$ so by proposition 3.2 (i) we get $\pi_{i+2,j+2}(sh(T))=\pi_{ij}(T)$,
$\forall\, T\in A_{ij}$.
$\qed$
 
\section{The universal monoidal category associated to a Popa
system}

We prove in this section the implication {\it (ii)}
$\Longrightarrow$ {\it (i)} in theorem 3.1. Let
$(A_{ij})_{0\leq i\leq j<\infty}$ be a Popa system, $H$ be a
finite dimensional Hilbert space and $\pi$ be a representation
of $(A_{ij})_{0\leq i\leq j<\infty}$ on $H$ such that
$\pi_{02}(e_2)$ is a rank one projection. We will construct a
pair $(A,u)$ satisfying {\it (w1-3)} such that
$(A_{ij})_{0\leq i\leq j<\infty}=End(u^{\otimes [i,j]})_{0\leq
i\leq j<\infty}$. This pair is far from being unique (see section 5); we will construct the ``universal'' one.

By proposition 3.1 we may assume that $\pi$ is normalised. Let $Q\in\l (H)$
be as in definition 3.2 and define linear maps
$i_{\alpha}$, $i_{\beta}$, $p_{\alpha}$,
$p_{\beta}$ as in proposition 3.2.

Let $\n^{*2}_{alt}$ be the subset of $\n^{*2}$ consisting of
alternating words, i.e. of words not containg $\alpha^2$ or
$\beta^2$. Thus $\n^{*2}_{alt}=\{ [i,j]\mid 0\leq
i\leq j\}$, and corollary 3.2 shows that to any
$x\in\n^{*2}_{alt}$ one can associate a $\c^*$-algebra
$A_x\subset \l (H^{\otimes x})$ by $A_x:=\pi_{i,j}(A_{i,j})$
for any $0\leq i\leq j$ such that $[i,j]=x$.

We define a monoidal subcategory ${\cal R}$ of the category of
finite dimensional complex vector spaces in the following way.
The objects of ${\cal R}$ are $\{ H^{\otimes x}\mid x\in
\n^{*2} \}$, and ${\cal R}$ is the smallest
monoidal category containing $i_{\alpha}$, 
$i_{\beta}$, $p_{\alpha}$,
$p_{\beta}$ and all the elements of all the $A_x$'s, for
$x\in \n^{*2}_{alt}$. Equivalently, the arrows of ${\cal R}$
are the linear combinations of (composable) compositions of
tensor products of maps of the form $i_{\alpha}$, 
$i_{\beta}$, $p_{\alpha}$,
$p_{\beta}$, or of the form $id_x:=id_{H^{\otimes x}}$
with $x\in\n^{*2}$, or of the form $T$ with $T\in A_x$ for
some $x\in \n^{*2}_{alt}$ (see also
\cite{w2},\cite{b0},\cite{b1} for this kind of constructions).

By definition of ${\cal R}$ we have $A_x\subset 
End_{\cal R}(H^{\otimes x})$ for any $x\in \n^{*2}_{alt}$. Theorem 3.1
will follow from the reconstruction results in \cite{w2} and from the following result.

\begin{prop}
$End_{\cal R}(H^{\otimes x})=A_x$ for any $x\in \n^{*2}_{alt}$.
\end{prop}

This will be proved after the ``computation'' of the arrows of
${\cal R}$ (lemma 4.1 and proposition 4.2). The following
special type of arrows will play an important role in this
computation.
\lvs
\noindent {\bf Definition 4.1} A block is a map of the form
$$id_a\otimes [(id_x\otimes
p_{\hat{\gamma}}^{\otimes n}\otimes id_z)T(id_y\otimes
i_\gamma^{\otimes m}\otimes id_t)]\otimes id_b$$
where $\gamma\in\{ \alpha ,\beta\}$, $m,n\geq 0$,
$a,b\in\n^{*2}$, $x,y,z,t\in\n^{*2}_{alt}$ and $T\in A_w$ for
some $w\in\n^{*2}_{alt}$ such that
$w=x(\gamma\hat{\gamma})^nz=y(\gamma\hat{\gamma})^mt$.
\lvs 

The above block will be denoted $B(a,b,x,y,z,t,\gamma
,n,m,w,T)$. It is an element of $Hom_{\cal R}(H^{\otimes
aytb}, H^{\otimes axzb})$. 

Such a block is said to be {\em connected} if
$a=b=e$;  {\em left degenerate} if $x=z=e$; {\em right
degenerate} if
$y=t=e$; {\em normalised} if $z=t=e$ and $mn=0$ (where
$e$ is the unit of $\n^{*2}$).

\begin{lemm}
Each arrow in ${\cal R}$ is a sum of compositions of blocks.
\end{lemm}

{\it Proof.} As the set of sums of compositions of blocks is
stable by linear sums and by compositions, we may consider the
category ${\cal R}^\prime$ whose objects are $\{ H^{\otimes
x}\mid x\in\n^{*2} \}$ and whose arrows are sums of
compositions of blocks. We have to show that ${\cal R}={\cal
R}^\prime$, and by definition of ${\cal R}$ and ${\cal
R}^\prime$ it is enough to show that ${\cal R}^\prime$ is a
monoidal category. But if $A\in\l (H^{\otimes x},H^{\otimes
y})$ and $B\in\l (H^{\otimes z},H^{\otimes
t})$ are blocks, then $A\otimes B$ is equal to $(id_y\otimes
B)(A\otimes id_z)$, which is a composition of blocks. $\qed$

\begin{lemm}
Each block is equal to a normalised block.
\end{lemm}

{\it Proof.} We will give an explicit method for normalising a
block. Let $B$ be as in the definition of blocks. By proposition 3.2 (i) the following formulas hold for any
$\gamma\in\{ \alpha ,\beta\}$:
$$i_\gamma\otimes id_\gamma =(i_\gamma p_{\hat{\gamma}}\otimes
id_\gamma )(id_\gamma\otimes i_{\hat{\gamma}})$$
$$p_\gamma\otimes id_{\hat{\gamma}} =
(id_{\hat{\gamma}}\otimes p_{\hat{\gamma}} )
(i_{\hat{\gamma}} p_\gamma\otimes id_{\hat{\gamma}} )$$
Also by proposition 3.2 (ii), $i_\gamma p_{\hat{\gamma}}$ and
$i_{\hat{\gamma}} p_\gamma$ are scalar multiples of
Jones projections, so by multiplying $T$ to the left and to
the right with suitable products of Jones projections, the
$i_\gamma$'s and $p_\gamma$'s can be moved to the right. That
is, $B$ is equal to a block $B^\prime$ having $z=t=e$.

If $mn=0$ we are done. Assume $m,n\geq 1$; then
$w=w^\prime\gamma\hat{\gamma}$ for some
$w^\prime\in\n^{*2}_{alt}$. Thus
$$B^\prime =id_a\otimes [(id_x\otimes
p_{\hat{\gamma}}^{\otimes n-1})T^\prime (id_y\otimes
i_\gamma^{\otimes m-1})]\otimes id_b$$
with $T^\prime=(id_{x(\gamma\hat{\gamma})^{n-1}}\otimes 
p_{\hat{\gamma}}
)T(id_{y(\gamma\hat{\gamma})^{n-1}}\otimes i_\gamma )$, which by corollary 3.1 is in $A_{w^\prime}$. By using this remark
$\mid m-n\mid$ times, it follows that $B^\prime$ is equal to a
normalised block. $\qed$

\begin{lemm}
A product of connected blocks is a connected block.
\end{lemm}

{\it Proof.} It is enough to do it for a product
$BB^\prime$ of two connected blocks. By lemma 4.2 we may
assume that $B$ and $B^\prime$ are connected and normalised.
Write 
$$B=(id_x\otimes p_{\hat{\gamma}}^{\otimes n})T(id_y\otimes
i_\gamma^{\otimes m}),\,\,\,  B^\prime =(id_{x^\prime} \otimes
p_{\hat{\gamma}^\prime}^{\otimes n^\prime} )T^\prime
(id_{y^\prime}
\otimes i_{\gamma^\prime}^{\otimes m^\prime})$$
with $nm=n^\prime m^\prime =0$. The composability of $B$ and
$B^\prime$ shows that $\gamma =\gamma^\prime$ and
$y=x^\prime$. Now each $i_\gamma p_{\hat{\gamma}}$ being a
scalar multiple of a Jones projection, by performing the
multiplication in the middle of
$BB^\prime$ we may assume that $m=0$ or $n^\prime =0$. Suppose
for instance that $m=0$; then $BB^\prime $ is equal to
$$(id_x\otimes p_{\hat{\gamma}}^{\otimes n+n^\prime})(T\otimes
id_{(\gamma\hat{\gamma})^{n^\prime}})T^\prime (id_{y^\prime}
\otimes i_{\gamma}^{\otimes m^\prime})$$
As $(T\otimes
id_{(\gamma\hat{\gamma})^{n^\prime}})T^\prime\in
A_{w^\prime}$, it follows that $BB^\prime$ is a connected block.
$\qed$

\begin{lemm}
Let $B$ be a right non-degenerate block and $B^\prime$ be a
left non-degenerate block. Write $B=id_a\otimes A\otimes id_b$
and $B^\prime =id_{a^\prime}\otimes A^\prime\otimes
id_{b^\prime}$ with $A,A^\prime$ connected blocks. If $B$ and
$B^\prime$ are composable, then $BB^\prime$ is a block or is
equal to a map of the form $id_f\otimes A\otimes id_g\otimes
A^\prime\otimes id_h$, or of the form $id_f\otimes
A^\prime\otimes id_g\otimes A\otimes id_h$ for some
$f,g,h\in\n^{*2}$.
\end{lemm}

{\it Proof.} Write $B=B(a,b,x,y,z,t,\gamma
,n,m,w,T)$ and $B^\prime=B(a^\prime$, $b^\prime$, $x^\prime$
, $y^\prime$, $z^\prime$,$t^\prime$,$\gamma^\prime$,$n^\prime$,
$m^\prime$,$w^\prime$,$T^\prime )$. As $B,B^{\prime}$ are
composable we get $aytb=a^\prime x^\prime z^\prime b^\prime$;
denote by $W$ this word. By the non-degeneracy assumptions,
both $yt$ and $x^\prime z^\prime$ are non-empty subwords of
$W$. There are two cases:

- either $yt$ and $x^\prime z^\prime$ are disjoint subwords of
$W$; in this case $BB^\prime$ is equal to $id_f\otimes
A\otimes id_g\otimes A^\prime\otimes id_h$ (if $yt$ is at the
right of $x^\prime z^\prime$) or to $id_f\otimes A^\prime\otimes
id_g\otimes A\otimes id_h$ (if $yt$ is at the left of $x^\prime
z^\prime$) for some $f,g,h\in\n^{*2}$.

- either $yt$ and $x^\prime z^\prime$ have at least one common
letter as subwords of $W$; we will prove in this case that 
$BB^\prime$ is a block. As
$w=y(\gamma\hat{\gamma})^mt$ and
$w^\prime
=x^\prime({\gamma^\prime}{\hat{\gamma}^\prime})^{n^\prime}t^\prime$
are alternating words, it follows that $yt$ and $x^\prime
z^\prime$ are alternating words, so their union in $W$ is an
alternating word, say $V$. Write
$W=a^{\prime\prime}Vb^{\prime\prime}$; then $a^{\prime\prime}$
is a subword of $a$, and if
$\varepsilon$ is such that $a=a^{\prime\prime}\varepsilon$ it
follows that $\varepsilon x\in\n^{*2}_{alt}$. Thus
$$B=id_{a^{\prime\prime}}\otimes [(id_{\varepsilon x}\otimes
p_{\hat{\gamma}}^{\otimes n}\otimes id_z)(id_\varepsilon\otimes
T)(id_y\otimes i_\gamma^{\otimes m}\otimes id_t)]\otimes id_b$$
with $id_\varepsilon\otimes T\in A_{\varepsilon w}$, so we may
suppose that $a=a^{\prime\prime}$. By the same argument we may
suppose $a^\prime =a^{\prime\prime}$, and also
$b=b^{\prime\prime}$ and $b^\prime =b^{\prime\prime}$. In this
way $BB^\prime$ becomes of the form
$id_{a^{\prime\prime}}\otimes K\otimes id_{b^{\prime\prime}}$,
with $K$ a product of connected blocks. By lemma 4.3 $K$
is a connected block, so $BB^\prime$ is a block. $\qed$

\begin{prop}
Every composition of blocks is of the form
$$C_1...C_nB_1...B_mA_1...A_p$$
where the $A_i$'s are left
degenerate blocks, the $C_i$'s are right degenerate blocks,
and the $B_i$'s are left and right non-degenerate blocks.
\end{prop}

{\it Proof.} Let us call ${\cal A}$ (resp. ${\cal C}$, ${\cal
B}$) the set of left degenerate (resp. right degenerate, left
and right non-degenerate) blocks. Lemma 4.4 shows that:

- if $A\in {\cal A}$ and $C\in {\cal C}$ then $AC$ is a block
or is of the form $C^\prime A^\prime$, with $A^\prime \in {\cal
A}$ and $C^\prime\in {\cal C}$.

- if $A\in {\cal A}$ and $B\in {\cal B}$ then $AB$ is a block
or is of the form $B^\prime A^\prime$, with $A^\prime \in {\cal
A}$ and $B^\prime\in {\cal R}$.

- if $B\in {\cal B}$ and $C\in {\cal C}$ then $BC$ is a block
or is of the form $C^\prime B^\prime$, with $B^\prime \in {\cal
B}$ and $C^\prime\in {\cal C}$.

Let $X$ be a composition of blocks, and
choose a decomposition of $X$ as product of minimal
number of blocks. By performing the above operations one may
write $X$ as in the statement. $\qed$

\lvs
\noindent {\bf Remark.} Lemma 4.1 and proposition 4.2 give the
structure of arrows in ${\cal R}$, as ``reduced words on
blocks''. Some better statements may be proved; however, there
is no simple writing for a product of type $A_1...A_p$ or
$C_1...C_n$.

\begin{lemm}
(i) Let $B\in Hom_{\cal R}(H^{\otimes u},H^{\otimes v})$ be a
right non-degenerate block. If $u\in \n^{*2}_{alt}$ then
$v\in\n^{*2}_{alt}$ and $B$ is equal to a connected block.

(ii) Let $B\in Hom_{\cal R}(H^{\otimes u},H^{\otimes v})$ be a
left non-degenerate block. If $v\in \n^{*2}_{alt}$ then
$u\in\n^{*2}_{alt}$ and $B$ is equal to a connected block.
\end{lemm}

{\it Proof.} We prove only (i). Write
$B=B(a,b,x,y,z,t,\gamma ,n,m,w,T)$. Then
$w=x(\gamma\hat{\gamma})^nz=y(\gamma\hat{\gamma})^mt$ (see the
definition of blocks),
$yt\neq e$ (by non-degeneracy) and $u=aytb$,
$v=axzb$ (as $B\in Hom_{\cal R}(H^{\otimes u},H^{\otimes v})$).
As
$u=aytb$ is an alternating word, it follows that:

- $a,b$ are alternating words.

- if $a$ ends with $\alpha$ (resp. $\beta$) then $yt$ begins
with $\beta$ (resp. $\alpha$).

- if $b$ begins with $\alpha$
(resp. $\beta$) then $yt$ ends with $\beta$ (resp. $\alpha$).

Also as $w=y(\gamma\hat{\gamma})^mt$, the alternating words $w$
and
$yt$ begin (resp. end) with the same letter. It follows that
$awb$ is an alternating word, and as $v$ is obtained from $awb$
by deleting some $\gamma\hat{\gamma}$'s, it is an alternating
word. Also
$awb\in \n^{*2}_{alt}$ shows that $B$ is equal to the connected
block 
$$(id_{ax}\otimes
p_{\hat{\gamma}}^{\otimes n}\otimes id_{zb})(id_a\otimes
T\otimes id_b)(id_{ay}\otimes i_\gamma^{\otimes m}\otimes
id_{tb})\,\,\, \qed$$

\begin{lemm}
If $x,y\in \n^{*2}_{alt}$ then the elements
$Hom_{\cal R}(H^{\otimes x},H^{\otimes y})$ are sums of
connected blocks.
\end{lemm}

{\it Proof.} Let $X\in Hom_{\cal R}(H^{\otimes
x},H^{\otimes y})$. By lemma 4.1 we may assume that $X$ is
a composition of blocks. Write 
$$X=C_1...C_nB_1...B_mA_1...A_p$$
as in proposition 4.2. As all $A_i$'s are left degenerate, we
may assume that they are right non-degenerate (a block which
is both left and right degenerate is a scalar). Thus lemma
4.5 (i) applies $p+m$ times and shows recursively that
$A_p,...,A_1,B_m,...,B_1$ are connected blocks. Also lemma 4.5
(ii) applies $n+m$ times and shows recursively that
$C_1,...,C_n,B_1,...,B_m$ are connected blocks. Thus $X$ is a
composition of connected blocks, and by lemma 4.3 it is a
connected block. $\qed$

\lvs
{\it Proof of proposition 4.1.} Lemma 4.5 shows that
every element $X\in End_{\cal R}(H^{\otimes x})$ is a sum of
connected blocks. By lemma 4.2, $X$ is a sum of normalised
connected blocks. But the only normalised connected blocks in 
$End_{\cal R}(H^{\otimes x})$ are the elements of $A_x$ (see the
definition of blocks), so we get $End_{\cal R}(H^{\otimes
x})\subset A_x$. The other inclusion is clear by construction
of ${\cal R}$. $\qed$

\lvs
{\em End of the proof of (ii) $\Longrightarrow$ (i) in
theorem 3.1.} We use here freely the terminology from
\cite{w2}. Let us view ${\cal R}$ as a concrete monoidal
(uncomplete) $\w^*$-category. The duality formulas in proposition 3.2 (i) show that the objects $H$ and $\Hb$ are conjugate
in ${\cal R}$. Thus theorem 1.3. in \cite{w2} applies and
shows that the ${\cal R}-$universal admissible pair $(A,u)$ is
a compact matrix pseudogroup (i.e. it satisfies the conditions
{\it (w1-3)} from the first section).  Let $w\in\l (\Hb
)\otimes A_s$ be the corepresentation corresponding to the
object $\Hb$. By the next lemma $w=\hat{u}$, so theorem 1.3 in \cite{w2} shows also that the space of
intertwiners $Hom(u^{\otimes x},u^{\otimes y})$ is
equal (as a subspace of $\l (H^{\otimes
x},H^{\otimes y})$) to the space $Hom_{\cal R}(H^{\otimes
x},H^{\otimes y})$ of arrows of ${\cal R}$, for any 
$x,y\in \n^{*2}$. In particular if $0\leq i\leq j$ then
$$End(u^{\otimes
[i,j]})=End_{\cal R}(H^{\otimes
[i,j]})=A_{[i,j]}=\pi_{i,j}(A_{i,j})$$
By the next lemma $Q_u=Q$, so the Jones
projections and the conditional
expectations of $(A_{ij})_{0\leq i\leq j<\infty}$
and of $End(v^{\otimes [i,j]})_{0\leq i\leq j <\infty}$ are
given by the same formulas (cf. proposition 3.2 (ii,iii) and the proof
of theorem 2.1). Thus the above identifications
$End(u^{\otimes [i,j]})\simeq A_{i,j}$ are Jones
projections-preserving and trace-preserving (the traces
being particular cases of conditional expectations). $\qed$

\begin{lemm}
$w=\hat{u}$ and $Q_u=Q$.
\end{lemm}

{\em Proof.} We first prove that $w=U\hat{u}U^*$ for some unitary $U$. We know from \cite{w2} that $w$ is equivalent
to $\hat{u}$. Let $F\in M_n(\c )$ be such that $w=F\hat{u}F^{-1}$. By polar decomposition, we may assume that $F$ is positive and we
have to prove that $w=\hat{u}$ in this case. As $w$ and $u$ are unitaries 
$$w=((id\otimes\kappa )w)^*=(F(id\otimes\kappa )(\hat{u})F^{-1})^*=F^{-1}((id\otimes\kappa )\hat{u})^*F=F^{-1}\hat{u}F$$
It follows that $F^2\in End(\hat{u})$, so $F\in End(\hat{u})$ and $w=F\hat{u}F^{-1}=\hat{u}$ as desired. 

As $i_{\alpha}\in Hom(1,u\otimes w)$ and
$i_{\beta}\in Hom(1,w\otimes u)$, it is easy to see that
$w=Q^t\ub (Q^{-1})^t$. On the other hand $w=U\hat{u}U^*$ and
$\ub =(Q_u^{-1})^t\hat{u}Q_u^t$, so $U^*Q^t(Q_u^{-1})^t\in
End(\hat{u})$. By polar decomposition both $U^*$ and
$Q^t(Q_u^{-1})^t$ are in
$End(\hat{u})$. As $U^*\in End(\hat{u})$, we get the first
assertion $w=U\hat{u}U^*=\hat{u}$. The other
relation $Q^t(Q_u^{-1})^t\in End(\hat{u})$ shows that 
$$[Q_u^t\ub (Q_u^{-1})^t]Q^t(Q_u^{-1})^t=Q^t(Q_u^{-1})^t[Q_u^t\ub (Q_u^{-1})^t]$$ 
It follows that $(QQ_u^{-1})^t\in End(\ub )$, so $QQ_u^{-1}\in
End(u)$. Thus the condition (iv) in lemma 1.2  is
satisfied. Also $Q$ is positive, so for proving $Q=Q_u$ it is
enough to show that $Tr(Q^{-2}.)=Tr(Q^2.)$ on $End(u)$. Indeed, proposition 3.2  shows that the expectation 
$\pi_{01}(A_{01})\rightarrow \pi_{00}(A_{00})=\c$ is
$Tr(Q^2.)$, and that the expectation 
$\pi_{01}(A_{01})\rightarrow \pi_{11}(A_{11})=\c$ is 
$Tr(Q^{-2}.)$. As both expectations coincide with the trace 
on $\pi_{01}(A_{01})$, and as $\pi_{01}(A_{01})=End(u)$, we get
$Tr(Q^{-2}.)=Tr(Q^2.)$ on
$End(u)$. $\qed$

\section{From Woronowicz algebras to Popa systems and back}

Let $X$ be the set of pairs (Woronowicz algebra,
corepresentation) and $Y$ be the set of pairs (Popa system,
normalised representation). In theorem 2.1 we constructed a
map $L:X\rightarrow Y$ by 
$$L(A,u)=(End(u^{\otimes
[i,j]})_{0\leq i\leq j<\infty},\pi )$$
where $\pi$ is the canonical representation (see example
3.1). The Theorem 3.1 says that $L$ is surjective. Moreover,
the proof of theorem 3.1 was as follows - to any
$a\in Y$ we have associated a certain category
${\cal R}$; then the universal ${\cal R}-$admissible pair
$R(a):=(A,u)$ was shown to satisfy $L(A,u)=a$. The map
$R:Y\rightarrow X$ being a section for $L$, it follows that
the composition
$$RL:X\rightarrow X$$
is a projection. We give in this section an
explicit description of $RL$. 

First of all we will give more precise definitions for
$X$ and $Y$. We begin with $X$. For simplicity we restrict
attention to the pairs $(A,u)$ satisfying the conditions {\it
(w1-3)} in the first section. We will use the following special
type of morphisms:
\lvs
\noindent {\bf Definition 5.1} If $(A,u)$ and $(B,v)$ satisfy {\it (w1-3)}, a strong
morphism from $(A,u)$ to $(B,v)$ is a $*$-algebra morphism
$f:A_s\rightarrow B_s$ such that $(id\otimes f)u=v$.
\lvs 

This definition has to be understood as follows. If $u\in
M_n(A_s)$ and $v\in M_m(B_s)$ with $m\neq n$ there is no strong
morphism between $(A,u)$ and $(B,v)$. If $m=n$ then there
exists {\em at most one strong morphism}, which has to send
$u_{ij}\mapsto v_{ij}$ for any $1\leq i,j\leq n$. Example: if
$G=<g_1,...,g_n>$ and $H=<h_1,...,h_m>$ are finitely generated
discrete groups, there exists a strong morphism
$$(\c^*(G), diag(u_{g_1},...,u_{g_n}))\rightarrow (\c^*(H),
diag(u_{h_1},...,u_{h_m}))$$
iff $m=n$ and there exists a group morphism $G\rightarrow H$
sending $g_i\mapsto h_i$ for every $i$. Note also that
$\c^* (G)$ and $\c^*_{red}(G)$ are strongly isomorphic; more
generally, given any $(A,u)$, if $A_p$ and $A_{red}$ are the
full and reduced version of $A$, then $(A,u)$, $(A_p,u)$,
$(A_{red},u)$ are strongly isomorphic (see
\cite{w1},\cite{bs}).

The following consequence of the uniqueness of strong
morphisms will be used several times: if
$f:(A,u)\rightarrow (B,v)$ and $g:(B,v)\rightarrow (A,u)$ are
strong morphisms, then $f$ and $g$ are both strong
isomorphisms.

We define $X$ to be the category of pairs $(A,u)$ satisfying
{\it (w1-3)}, with the strong morphisms. $Y$ is by definition
the set of quadruples consisting of a Popa
system $(A_{ij})_{0\leq i\leq j<\infty}$, a Hilbert space $H$,
a positive operator $Q\in \l (H)$, and a representation $\pi$
of $(A_{ij})_{0\leq i\leq j<\infty}$ on $H$ as in definition 3.2, with the obvious notion of equality for such
quadruples. Let us first give the abstract description of the
map $RL$.

\begin{prop}
To any pair $(A,u)$ satisfying {\it (w1-3)} we associate a
category ${\cal C} (A,u)$ in the following way:

- the objects of ${\cal C} (A,u)$ are the pairs $(B,v)$
satisfying {\it (w1-3)} and such that $L(B,v)=L(A,u)$. 

- the arrows of ${\cal C} (A,u)$ are the strong morphisms.

Then $RL(A,u)$ is the (unique) universally repelling object of
${\cal C} (A,u)$.
\end{prop}

{\it Proof.} Consider the Popa system $End(u^{\otimes
[i,j]})_{0\leq i\leq j<\infty}$ together with its
canonical normalised representation $\pi$. Let ${\cal R}$ be
the category defined in the fourth section, so that $RL(A,u)$
is the universal ${\cal R}-$admissible pair. By definition of
${\cal R}$ and of ${\cal C} (A,u)$ we see that every element of
${\cal C} (A,u)$ is an ${\cal R}-$admissible pair, so the
assertion is just a translation of the universal property of the
universal admissible pair (the unicity up to strong
isomorphism is clear from the unicity of strong morphisms).
$\qed$

\lvs
We will need the following results on free products of
discrete quantum groups \cite{wa}. If $A$ and $B$ are
Woronowicz algebras, so is their free product $A*B$ (=
coproduct in the category of unital
$\c^*$-algebras). $Irr(A*B)-\{ 1\}$ is then the set of
alternating products of elements of
$Irr(A)-\{ 1\}$ with elements of $Irr(B)-\{ 1\}$. The Haar 
measure of $A*B$ is the free product $h*k$ of Haar measures $h$
of $A$ and $k$ of $B$. If $*_{red}$  denotes the reduced free
product with respect to $h*k$, then $A*_{red}B$ is also a
Woronowicz algebra (which is equal to the reduced version of
the Woronowicz algebra $A*B$). The
$\c^*$-algebras
$A_{red}$ and
$B_{red}$ are embedded in $A*_{red}B$, and are free in the
sense of \cite{vdn} with respect to $h*k$.

\lvs 
Consider the Woronowicz algebra $\c^*(\z )$ and let $z$ be the
unitary of $\c^*(\z )$ corresponding to the generator $1$ of
$\z$; it is a one-dimensional corepresentation of $\c^*(\z )$.
Recall that we have a canonical isomorphism $\c^*(\z )\simeq C(\t )$, which maps $z$ to the
function $\t\ni x\mapsto x\in\c $. 

\begin{theo}
To any pair $(A,u)$ satisfying {\it (w1-3)} we associate a
pair $(\tilde{A},\tilde{u})$ satisfying {\it (w1-3)} in the
following way: $\tilde{A}$ is the
$\c^*$-subalgebra of $\c^*(\z )*A$ generated by the entries of
the matrix $\tilde{u}:=zu$. Then $(\tilde{A},\tilde{u})$ is
(strongly isomorphic to) the universally repelling object
$RL(A,u)$ of ${\cal C} (A,u)$.
\end{theo}

We begin with a few remarks on the operation $(A,u)\mapsto
(\tilde{A},\tilde{u})$. First, as $\tilde{u}$ is a
corepresentation of $\c^*(\z )*A$, the comultiplication and the
antipode of $\tilde{A}$  are the restrictions of the ones of
$\c^*(\z )*A$. Note however that the $\c^*$-algebra $\tilde{A}$
depends on both $u$ and $A$. The reduced version
$\tilde{A}_{red}$ is the $\c^*$-subalgebra of $\c^*(\z
)*_{red}A$ generated by the entries of the  matrix $zu$. Note
that the full version $\tilde{A}_p$ may be different from
$\tilde{A}$. Of course,
$(\tilde{A}_p,\tilde{u})$, $(\tilde{A},\tilde{u})$ and
$(\tilde{A}_{red},\tilde{u})$ are strongly isomorphic.

\begin{lemm}
(i) $(A,u)\mapsto (\tilde{A},\tilde{u})$ is functorial.

(ii) There exists a strong morphism
$(\tilde{A},\tilde{u})\rightarrow (A,u)$.

(iii) $(\,\,\tilde{\!\!\tilde{A}},
\tilde{\tilde{u}})$ and $(\tilde{A},\tilde{u})$ are strongly
isomorphic.
\end{lemm}

{\em Proof.} As there is at most one strong morphism between
two objects, point (i) states precisely that if $(A,u)$
and $(B,v)$ are such that there exists a strong morphism
$f:(A,u)\rightarrow (B,v)$, then there exists a strong
morphism $(\tilde{A},\tilde{u})\rightarrow
(\tilde{B},\tilde{v})$. Such a strong morphism $\tilde{f}$ may
be constructed in the following way
: $\tilde{f}$ is the restriction to $\tilde{A}_s$ of the
$\c^*$-morphism $\psi :\c^*(\z )*A\rightarrow \c^*(\z )*B$
defined by
$\psi\mid_{\c^*(\z )}=id$ and $\psi\mid_A=f$. 

A strong morphism as in (ii) could be
constructed as  the restriction to $\tilde{A}_s$ of the
$\c^*$-morphism $\psi :\c^*(\z )*A\rightarrow A$  defined by
$\psi\mid_{\c^*(\z )}=\varepsilon\mid_{\c^*(\z )}$ (the
counit of $\c^*(\z )$) and $\psi\mid_A=id_A$. 

For (iii) note that
$\,\,\tilde{\!\!\tilde{A}}_s$ is  the $*$-subalgebra of
$\c^*(\z )*\c^*(\z )*A$ generated by the entries of
$\tilde{\tilde{u}} =z^\prime zu$, so a strong morphism 
$(\tilde{A},\tilde{u})\rightarrow (\,\,\tilde{\!\!\tilde{A}}, \tilde{\tilde{u}})$ could be constructed as the restriction to 
$\tilde{A}_s$ of the $\c^*$-morphism $\phi :\c^*(\z )*A\rightarrow \c^*(\z )*\c^*(\z )*A$ defined by 
$\phi\mid_{\c^*(\z )}=z^\prime z\varepsilon\mid_{\c^*(\z
)}$ and $\phi\mid_A=id_A$. Now (iii) follows from (ii)
and from the unicity of the strong morphisms. $\qed$

\lvs
We have to prove that the functor $(A,u)\mapsto
(\tilde{A},\tilde{u})$ ``keeps fixed the associated Popa
system'' and ``destroys the rest of the structure''. The
first assertion is easy and is the next lemma; the other
assertion will follow from proposition 5.1 and from the
isomorphism criterion in proposition 5.2. 

\begin{lemm}
$Q_{\tilde{u}}=Q_u$, $\widehat{\tilde{u}}=\hat{u}z^*$, and
$L(\tilde{A},\tilde{u})=L(A,u)$.
\end{lemm}

{\em Proof.} The first two assertions are clear from lemma
1.2. Now using $\tilde{u}=zu$, $\widehat{\tilde{u}}=\hat{u}z^*$
and $zz^*=z^*z=1$ we find that for every $x\in\n^{*2}_{alt}$,
$$\tilde{u}^{\otimes x}=z^mu^{\otimes x}z^n$$
for some $m,n\in \{ -1,0,1\}$. Thus $End(\tilde{u}^{\otimes
x})=End(u^{\otimes x})$ for any alternating word $x$, and this
proves $L(\tilde{A},\tilde{u})=L(A,u)$. $\qed$

\begin{lemm}
Let $f:(A,u)\rightarrow (B,v)$ be a strong
morphism. If 
$$dim(Hom(u^{\otimes
x},u^{\otimes y}))= dim(Hom(v^{\otimes x},v^{\otimes y}))$$
for any $x,y\in\n^{*2}$, then $f$ is a strong 
isomorphism.
\end{lemm}

{\it Proof.} Let $f_*$ be the map induced by $f$ at the
level of classes of corepresentations. Let $y$ be an arbitrary
element of $\n^{*2}$. Write
$u^{\otimes y}=\sum m_ir_i$, with $m_i>0$ and $r_i$
irreducible and disjoint. Then 
$$v^{\otimes y}=(f_*u)^{\otimes y}=f_*(u^{\otimes y})=\sum
m_if_*(r_i)$$
By \cite{w1} we have $dim(End(u^{\otimes y}))=\sum m_i^2$, and
$dim(End(v^{\otimes y}))\geq\sum m_i^2$, with equality iff the
$f_*(r_i)$'s are irreducible and disjoint. As
$dim(End(u^{\otimes y}))=dim(End(v^{\otimes y}))$, the
$f_*(r_i)$'s have to be irreducible. 

Also from
$dim(Hom(1,u^{\otimes y}))= dim(Hom(1,v^{\otimes y}))$ we get
that $r_i\neq 1$ implies $f_*(r_i)\neq 1$. As every irreducible
corepresentation of $A$ appears in some tensor product of the
form $u^{\otimes y}$ (see \cite{w1}) we find that $f_*$ maps
$Irr(A)-\{ 1\}$ into $Irr(B)-\{ 1\}$. 

Now let $a\in A_s$ be arbitrary. Write $a=\lambda 1+\sum a_i$,
with $\lambda\in\c$ and with $a_i$= coefficients of
irreducible non-trivial corepresentations of $A$ (cf.
\cite{w1}). Then
$f(a)=\lambda 1+\sum f(a_i)$, and by the above result the
$f(a_i)$'s are coefficients of
irreducible non-trivial corepresentations of $B$. Now if
$h_A,h_B$ are the Haar measures of $A$ and $B$, it follows that
$h_B(f(a))=\lambda =h_A(a)$. Thus the equality $h_Bf=h_A$ holds
on $A_s$, and by positivity we get that the restriction of $f$
to $A_s$ is injective. $\qed$

\begin{prop}
Let $f:(A,u)\rightarrow (B,v)$ be a strong
morphism. If 
$$dim(End(u^{\otimes
x}))=dim(End(v^{\otimes x}))$$
for any $x\in\n^{*2}_{alt}$, then
$(\tilde{A},\tilde{u})$ and $(\tilde{B},\tilde{v})$ are
strongly isomorphic.
\end{prop}

{\it Proof.} By functoriality (lemma 5.1. (i)), the conclusion
will follow from lemma 5.3 once we will prove that 
$$dim(Hom(\tilde{u}^{\otimes x},\tilde{u}^{\otimes
y}))=dim(Hom(\tilde{v}^{\otimes x},\tilde{v}^{\otimes y}))$$
for any $x,y\in\n^{*2}$. In turn, this is equivalent to the
following statement: 

\lvs
{\em If $(A,u)$ satisfies {\it (w1-3)}
then the numbers $dim(Hom(\tilde{u}^{\otimes x},\tilde{u}^{\otimes
y}))$ (with $x,y\in\n^{*2}$) depend only on the numbers
$dim(End(u^{\otimes x}))$ (with $x\in\n^{*2}_{alt}$).}

\lvs
We will use the following key remark (see also \cite{b1}). If
$(B,v)$ satisfies {\it (w1-3)} then by the orthogonality
formulas for characters in
\cite{w1} we have
$$dim(Hom(v^{\otimes x},v^{\otimes y}))=h\chi (v^{\otimes
x\overline{y}})=h(\chi (v)^{x\overline{y}})\,\,\,\,\,
(\star )$$ where $h\in B^*$ is the Haar measure and $\chi$ is
the character of corepresentations (and where $\chi
(v)^z$ denotes the image of $z\in\n^{*2}$ by
the morphism of involutive monoids $\n^{*2}\rightarrow (B,\cdot
,*)$ defined by $\alpha\mapsto \chi (v)$). Thus the numbers
$(\star )$ are the $*$-moments of the non commutative random
variable $\chi (v)\in (B,h)$ in the sense of \cite{vdn}.

By \cite{wa}, $z$ and $\chi (u)$ are $*$-free in $\c^*(\z
)*A$ with respect to its Haar measure, so $(\chi
(\tilde{u}),\chi (\tilde{u})^*)=(z\chi (u),\chi (u)^*z^*)$ is
an $R$-diagonal pair in the sense of \cite{ns}. By
\cite{ns}, formula (1.12) the determining series of the
$R$-transform of $(\chi (\tilde{u}),\chi (\tilde{u})^*)$ depends
only on the $R$-transform of $\chi (u)\chi (u)^*$. This
means that the $*$-moments of $\chi (\tilde{u})$
depend only on the moments of $\chi (u)\chi (u)^*$. But the
above remark shows that the $*$-moments of $\chi (\tilde{u})$
are the numbers $dim(Hom(\tilde{u}^{\otimes
x},\tilde{u}^{\otimes y}))$ (with $x,y\in\n^{*2}$), and that
the moments of $\chi (u)\chi (u)^*$ are the numbers
$dim(End(u^{\otimes x}))$ with
$x\in\n^{*2}_{alt}$. $\qed$

\lvs
{\it Proof of theorem 5.1.} Let $(B,v)$ be the universally
repelling object of ${\cal C} (A,u)$. By its universal
property we get a strong morphism $(B,v)\rightarrow (A,u)$, and
as $L(B,v)=L(A,u)$, proposition 5.2 applies and shows
that $(\tilde{A},\tilde{u})$ and $(\tilde{B},\tilde{v})$ are
strongly isomorphic. 

Lemma 5.2 shows that $(\tilde{B},\tilde{v})\in {\cal C}
(A,u)$, so by the universal property of $(B,v)$ we get a strong
morphism $(B,v)\rightarrow (\tilde{B},\tilde{v})$. Together
with lemma 5.1 (ii) and with the unicity of strong
morphisms, this shows that $(B,v)$ and $(\tilde{B},\tilde{v})$
are strongly isomorphic. $\qed$ 

\lvs
The operation $(A,u)\mapsto (\tilde{A},\tilde{u})$ seems to be
interesting, and we end this section with a complete
computation for discrete groups, and with a few remarks, to be
proved somewhere else. 

\begin{prop}
If $\Gamma$ is a discrete group generated by $g_1,...,g_n$ then 
$$(\tilde{\c}^*(\Gamma ),
\tilde{d}iag(u_{g_1},...,u_{g_n}))\simeq
(\c^*(\tilde{\Gamma}),diag(u_{zg_1},...,u_{zg_n}))$$ 
where $\tilde{\Gamma}$ is the subgroup of $\z *\Gamma$ generated
by $zg_1,...,zg_n$ and $z$ is the generator of $\z$. The
group $\tilde{\Gamma}$ could be computed as follows - if $H$ is
the subgroup of $\Gamma$ generated by 
$\{ g_i^{-1}g_j\mid i,j=1...n\}$ then:

(i) if $H=\Gamma$ then $\tilde{\Gamma}=\z *\Gamma$.

(ii) if $H\neq\Gamma$ then there exists an isomorphism $\tilde{\Gamma}\rightarrow \z *H$ sending $zg_1\mapsto z$ and 
$zg_i\mapsto g_1^{-1}g_i$ for $i=2,...,n$.
\end{prop}

{\em Proof.} As
$\tilde{\Gamma}=<zg_1,(zg_1)(g_1^{-1}g_2),...,(zg_1)(g_1^{-1}g_n)>=<zg_1,H>$
there are two cases:

(i) if $H=\Gamma$ then $g_1\in H$, so $z\in
\tilde{\Gamma}$, so $\tilde{\Gamma}=\z *\Gamma$.

(ii) if $H\neq \Gamma$ then $g_1\notin H$. Let
$Y$ be the subgroup generated in $\z *\Gamma$ by $zg_1$. We
have to prove that $Y$ and $H$  are free in $\z *\Gamma$.
Suppose that there exist $y_i\in Y-\{ 1\}$ and $h_i\in H-\{
1\}$ such that
$y_1h_1y_2h_2...=1$.  Choose such a product $P$ having a
minimal number of $z$'s and $z^{-1}$'s in its decomposition.
It's easy to see that $P$ is  a product of $z,z^{-1}$'s
alternating with terms of the form
$g_1,g_1^{-1},h_i,g_1h_i,h_ig_1^{-1}$ or $g_1h_ig_1^{-1}$. As
$P=1$ and 
$z$ is free from $\Gamma$, at least one of these terms has to be equal to $1$. But the only terms which can be equal to $1$ 
are the ones of the form $g_1h_ig_1^{-1}$, which can only appear between a $z$ and a $z^{-1}$. Now by deleting 
$zg_1h_ig_1^{-1}z^{-1}$ from $P$ we get another product equal to
$1$ which has less $z$'s and $z^{-1}$'s, 
contradiction. $\qed$

\lvs
\noindent {\bf Remarks.} For $A=C(G)$ with $G$ compact non-abelian the
Woronowicz algebra $\tilde{A}$ is harder to describe. Note
however that by \cite{b1} we have $\tilde{C}(\s\u (2))=
A_u(I_2)$; in fact this result, as well as its
generalisation $\tilde{A}_o(F)=A_u(F)$ may be
deduced from theorem 5.1. 

For any pair $(A,u)$ satisfying {\it (w1-3)} let $A^-$ be
the $\c^*$-subalgebra of $A$ generated by the the entries of
the matrix $u^- :=\hat{u}\otimes u$. Then $(A^-
,u^- )$ satisfies {\it (w1-3)} and it is easy to see
that $(\tilde{A}^-,\tilde{u}^- )=(A^-
,u^- )$. It is interesting to note that in proposition
5.3, $\tilde{\Gamma}$ is isomorphic to $\z *H$ in all cases.
This could be interpreted in the
following way: if $(A,u)$ is of the form $(\c^*(\Gamma
),diag(u_{g_1},...,u_{g_n}))$ then $\tilde{A}\simeq \c^*(\z
)*A^-$. This kind of free product decomposition of
$\tilde{A}$ is not available in general, but is valid for instance
if the trivial representation is contained in an alternating
product of the form $u\otimes \hat{u}\otimes ...\otimes
\hat{u}$ (the proof is similar to the one of proposition 5.3 (i)). 

It is also possible in some of the remaining cases to obtain
free product decompositions $vN(\tilde{A})\simeq vN(\z
)*vN(A^-)$ at the level of von Neumann algebras (see theorem 6 in \cite{b1}; the argument in its proof works for any
$(A,u)$ such that the Haar measure of $A$ is a trace and such
that the polar part of some coefficient of $u$ is unitary). 

The simplicity results in \cite{b1} may also be
extended - if the trivial representation is not contained in
any alternating product of the form $u\otimes \hat{u}\otimes
...\otimes \hat{u}$, one can show (first by studying the
operation $Irr(A)\mapsto Irr(\tilde{A})$, then by applying proposition 8 in \cite{b1}) that the reduced version of
$\tilde{A}$ is simple, with at most one trace.
\section{Amenability}

In this section we prove a Kesten type result for the (discrete
quantum groups represented by) Woronowicz algebras. This will
give a characterisation of the amenable Popa systems of the form
$End(u^{\otimes [i,j]})_{0\leq i\leq j <\infty}$. 
\lvs
\noindent {\bf Definition 6.1} Let $A$ be a Woronowicz algebra and let $h\in A^*$ be
its Haar measure. Define $A_{red}=A/\{ x\in A\mid h(x^*x)=0\}$
and $A_p=$ envelopping $\c^*$-algebra of 
$A_s$. Then $A_p$ and $A_{red}$ are Woronowicz algebras, called
the full and the reduced version of $A$ \cite{w1},\cite{bs}.
$A$ is said to be amenable (as a Woronowicz algebra) if the
canonical surjection $A_p\rightarrow A_{red}$ is an
isomorphism.
\lvs

\noindent {\bf Examples.} If $\Gamma$ is a discrete group then the
Woronowicz algebra $\c^*(\Gamma )$ is amenable if and only if
$\Gamma$ is amenable. If $G$ is a compact group then
the Woronowicz algebra $C(G)$ is clearly amenable (while this
is not related to the fact that $G$ is amenable). See
\cite{bs},\cite{bl} for the analogue of the Pontriagyn duality
and for the two dual notions of amenability. 

We restrict the attention to the Woronowicz algebras
satisfying {\it (w1-3)} for some corepresentation $u$ (these
are the ones which represent the discrete quantum groups ``of
finite type''). We recall that the character of a finite dimensional corepresentation
$r\in \l (V)\otimes A_s$ is $\chi (r)=(Tr\otimes id)r$.

\begin{theo} (G. Skandalis) 
Assume that $(A,u)$ satisfies {\it (w1-3)} with $u\in
M_n(A)$, and let $X\subset \r$ be the spectrum of $Re(\chi
(u))\in A_{red}$.

(i) $X\subset [-n,n]$.

(ii) The Woronowicz algebra $A$ is amenable if and only if $n\in
X$.
\end{theo}

{\it Proof.}  {\it (i)} is clear, as all entries of the unitary
$u$ are of norm $\leq 1$. 

{\it (ii,$\Longrightarrow )$} The counit of $A_s$ extends to
the full version as a $\c^*$-morphism 
$A_p\rightarrow\c$ which sends $n-Re(\chi (u))\mapsto 0$. 
Thus $n-Re(\chi (u))$ is not invertible in $A_p$, nor in
$A_{red}=A_p$ if $A$ is amenable. 
\lvs
{\it (ii,$\Longleftarrow )$} Embed
isometrically $A_{red}\hookrightarrow B(l^2(A))$ by the left
regular representation. The matrix $u$ being unitary, its
coefficients have norms $\leq 1$. In particular for all 
$i$ the operator $\alpha_i:=1-Re(u_{ii})$ is positive.

Assume that $n$ is in the spectrum of $Re(\chi (u))$, i.e. that the positive operator 
$\sum_i \alpha_i$ is not invertible. Then there 
exists a sequence $\xi_k$ of norm one vectors in $l^2(A)$ such that $<(\sum_i \alpha_i)\xi_k,\xi_k>\rightarrow 0$. 
As all $<\alpha_i \xi_k,\xi_k>$'s are positive, we get $<\alpha_i\xi_k,\xi_k>\rightarrow 0$ for all $i$ and 
by using the equality 
$$\no (1-x)\xi\no^2=2<(1-Re(x))\xi ,\xi >+\no x\xi\no^2 -1$$
valid for every $x$ and norm one vector $\xi$ we get $\no
u_{ii}\xi_k-\xi_k\no\rightarrow 0$ for all $i$. 

Now let $M_n(A_{red})$ act on $\c^n\otimes l^2(A)$. As $u$ is unitary $\no u(e_j\otimes
\xi_k)\no =1$ for all $j$ (where $\{ e_1,...,e_n\}$ is the canonical basis of $\c^n$). Thus 
$\sum_i \no u_{ij}\xi_k\no^2=1$, and as $\no u_{ii}\xi_k-\xi_k\no\rightarrow 0$ we get 
$\no u_{ij}\xi_k\no\rightarrow 0$ if $i\neq j$. 

These results could be summarized as $\no u_{ij}\xi_k-\delta_{ij}\xi_k\no\rightarrow 0$ for 
all $i,j$. Now let $\varepsilon :A_p\rightarrow\c$ be the
counit, $\pi :A_p\rightarrow A_{red}$ be the canonical
projection and  consider the set $M$ of elements of $A_p$
satisfying 
$$\no\pi (x)\xi_k-\varepsilon (x)\xi_k\no\rightarrow 0$$
$M$ is clearly a closed $*$-subalgebra of $A_p$, and contains all the coefficients of $u$, so $M=A_p$ and proposition 5.5 in \cite{bl} shows that $A$ is amenable.
$\qed$

\begin{theo}
If $v$ is a corepresentation of a Woronowicz algebra $A$ and if
$B$ is the $\c^*$-subalgebra of $A$ generated by the
coefficients of $v\otimes \hat{v}$ then the following are
equivalent:

(i) $End(v^{\otimes [i,j]})_{0\leq i\leq j <\infty}$ is
amenable.

(ii) The Woronowicz algebra $B$
is amenable and its Haar measure is a trace.
\end{theo}

{\em Proof.} Recall that a Popa $\lambda$-system is amenable if
the norm of its principal graph is $\lambda^{-1/2}$. Let
$\Gamma$ be the principal graph of $End(v^{\otimes
[i,j]})_{0\leq i\leq j <\infty}$ and $\lambda :=d(v)^{-2}$.
Let also $u:=v\otimes \hat{v}$ and consider the lattice
$$End(u^{\otimes
[i,j]})_{0\leq i\leq j <\infty} = End(v^{\otimes
[i,j]})_{0\leq i\leq j <\infty ,i+j\, even}$$
and denote by $\Gamma^\prime$ its principal
graph. Then $\no \Gamma^\prime\no
=\no \Gamma\no^2$, and $d(u)=d(v)^2$ by lemma 1.3, so the
amenability of $End(v^{\otimes [i,j]})_{0\leq i\leq j <\infty}$
is equivalent to the equality
$$\no \Gamma^\prime \no =d(u)$$
As $1\in u=\hat{u}$ (as equivalence classes of
corepresentations) we see that every element of
$Irr(B)$ is a subcorepresentation of $u^{\otimes n}$,
$\forall\, n\geq N$, for some $N\in\n$. By identifying the
vertices of $\Gamma^\prime$ with minimal
central projections in the algebras $End(u^{\otimes n})$, thus
with irreducible subcorepresentations of $u^{\otimes n}$, we see
that both sets of even and odd vertices are equal to
$Irr(B)$. Thus the norm of $\Gamma^\prime$ is equal to the
norm of the matrix 
$$M:=(dim(Hom(r\otimes u,p)))_{r,p\in Irr(B)}$$
Now $dim(Hom(r\otimes u,p))=<\chi (u)\chi (r),\chi (p)>$,
where $<,>$ is the scalar product associated to the Haar
measure on the algebra of characters $B_{central}$ (see
\cite{w1}). Thus $M$ is the multiplication by $\chi (u)$ in
$B_{central}$, so the norm of
$M$ is equal to the norm of $\chi (u)$ in $B_{red}$. After
all these identifications, we have the following equalities and
inequalities (where the inequalities are clear):
$$\no \Gamma^\prime\no =\no M\no =\no \chi
(u)\no_{B_{red}}\leq dim(u)\leq d(u)$$
Now $dim(u)=d(u)$ iff the Haar measure of $B$ is a
trace (see example 1.1), and by theorem 6.1 $\no\chi (u)\no_{B_{red}}=dim(u)$ iff $B$ is
amenable (note that $\chi (u)=\chi (v)\chi (v)^*$ is
positive). $\qed$

\begin{coro}
If $End(v^{\otimes [i,j]})_{0\leq i\leq j
<\infty}$ is amenable then its index is
the square of an integer. $\qed$
\end{coro}

It is somehow clear from the proof of theorem 6.2 that the
notion of amenability of a Woronowicz algebra $A$ depends only
on its corepresentation theory. This remark will give us the
co-amenability of the $q$-deformations. More precisely,
this will follow from the following consequence of theorem 6.1,
which seems to have its own interest: it gives for instance
information (a lower bound for dimensions) on the fiber
functors on the monoidal category of corepresentations of an
amenable Woronowicz algebra (see also examples 1.3., 1.4,
and the concluding remarks in the introduction).

We recall that if $A$ is a Woronowicz algebra, the fusion semiring
$R^+(A)$ is the semiring whose elements are equivalence classes of finite
dimensional corepresentations of $A$, and whose operations are
the sum and tensor product of (classes of) corepresentations.

\begin{prop}
Assume that $(A,u)$ and $(B,v)$ satisfy {\it (w1-3)} and that
there exists an isomorphism of semirings $f:R^+(A)\rightarrow
R^+(B)$ sending $u\mapsto v$. If the Woronowicz algebra $A$ is
amenable then $dim(f(r))\geq dim(r)$, $\forall\, r\in R^+(A)$,
and the following conditions are equivalent.

(i) $dim(f(r))=dim(r)$, $\forall\, r\in R^+(A)$.

(ii) $dim(v)=dim(u)$.

(iii) The Woronowicz algebra $B$ is amenable.
\end{prop}

{\em Proof.} {\it Step I.} We prove that $Spec(Re(\chi
(u)))=Spec(Re(\chi (v)))$. By \cite{w1} $R^+(A)$ and $R^+(B)$
are (as additive monoids) the free monoids on $Irr(A)$ and
$Irr(B)$ respectively, and it follows that $f$ maps $Irr(A)$
onto $Irr(B)$. Let us prove now that $f$ is an isomorphism of
involutive semirings, i.e. that $f(\hat{r})=\widehat{f(r)}$ for
any $r\in R^+(A)$. For $r\in Irr(A)$ this is clear, as
$\hat{r}$ may be characterised as being the unique $p\in
Irr(A)$ such that there exists $q\in R^+(A)$ with $r\otimes
p=q+1$. The general case follows by complete reducibility.

Now by functional calculus the spectrum of $Re(\chi (u))\in
A_{red}$ is the support of its spectral measure $\mu$ with
respect to the Haar measure $h$ of $A_{red}$, which in turn is
uniquely determined by its moments 
$$\mu (X^k)=h[Re(\chi
(u))^k]=2^{-k}dim(Hom(1,(u+\hat{u})^{\otimes k}))$$
hence by the isomorphism class of the pointed involutive
semiring $(R^+(A),u)$. The equality of spectra follows.

{\it Step II.} We prove that $dim(f(r))\geq dim(r)$, $\forall\,
r\in R^+(A)$. Firstly, by theorem 6.1 and by step {\it
I.} we get $dim(v)\geq dim(u)$. 

For any corepresentation $r$ of $A$ let
$\c^*(r)$ be the $\c^*$-subalgebra of $A$ generated by the
coefficients of $r$.  Then (by fixing a basis making $r$
unitary) $(\c^*(r),r)$ satisfies {\it (w1-3)}. As $A$ is
amenable, $\c^*(r)$ is also amenable (see \cite{bs},\cite{bl},
in  fact this could be deduced also from theorem 6.1). The
inequality $dim(f(r))\geq dim(r)$ follows by replacing $(A,u)$
with $(\c^*(r),r)$. 

{\it Step III.} The equivalences are clear: {\it (i)
$\Longrightarrow$ (ii)} is trivial, {\it (ii) $\Longrightarrow$
(iii)} follows from theorem 6.1 and from step {\it I.},
and if {\it (iii)} is satisfied, then by interchanging $(A,u)$
and $(B,v)$ we get the reverse inequalities in {\it (i)}. $\qed$

\begin{coro}
If $G$ is one of the groups $SU(n)$,
$Sp(2n)$ or $Spin(n)$ and if $q>0$ then the Woronowicz algebra
$C(G_q)$ is amenable.
\end{coro}

{\it Proof.} By \cite{r1},\cite{ro} we have an isomorphism
$R^+(C(G_q))\simeq R^+(C(G))$ which preserves the dimensions.
The assertion follows from {\it (iii)
$\Longleftrightarrow$ (i)} in proposition 6.1 and from the fact
that the Woronowicz algebra $C(G)$ is amenable (see definition 6.1). $\qed$

\lvs
\noindent {\bf Note.} For $SU(n)$ this was done by Nagy in \cite{n}.

\bigskip\noindent Teodor Banica

\noindent Institut de Math\'ematiques de Luminy,
Marseille, and

\noindent Institut de Math\'ematiques de Jussieu, case 191,
Universit\'e Paris 6, 4 place Jussieu, 75005 Paris, France

\noindent {\em E-mail:} banica@math.jussieu.fr
\end{document}